\documentclass{amsart}

\usepackage[letterpaper,text={6.5in,9in},centering]{geometry}

\usepackage[numbers, sort]{natbib}

\usepackage{amsmath}
\usepackage{amssymb}
\usepackage{amsmath}
\usepackage{amsfonts}
\usepackage{enumerate}
\usepackage{amsthm}
\usepackage{setspace}
\usepackage{multirow}
\usepackage{mathdots}
\usepackage{mathtools}
\usepackage{graphicx}
\usepackage{float}
\usepackage{subcaption}

\graphicspath{{./Figures/}}

\author{Xuefeng Shen}
\author{Melvin Leok}
\title{Geometric Symmetry Reduction of the Unobservable Subspace for Kalman Filtering}

\begin{document}

\begin{abstract}
In this article, we consider the implications of unobservable subspaces in the construction of a Kalman filter. In particular, we consider dynamical systems which are invariant with respect to a group action, and which are therefore unobservable in the group direction. We obtain reduced propagation and measurement equations that are invariant with respect to the group action, and we decompose the state space into unobservable and observable parts. Based on the decomposition, we propose a reduced Bayesian inference method, which exhibits superior accuracy for orientation and position estimation, and that is more robust to large measurement noise.
\end{abstract}

\maketitle

\section{Introduction}
Kalman filters are widely used in practical state estimation problems, including SLAM (Simultaneous Localization and Mapping). This is primarily due to its conceptual simplicity and the low computational complexity compared to optimization based methods. Various types of sensors are used for localization tasks, including GPS, odometry, inertial measurement units (IMU), cameras, Light Detection and Ranging (LIDAR), and so on, and the specific choice of sensors used depends on the application. Kalman filters have been successfully applied to Visual Odometry (VO)~\cite{monoslam2007}, Visual Inertial Odometry (VIO)~\cite{msckf2007}, and more recently, to autonomous driving cars, where information from multiple sensors (GPS, odometry, IMU, and cameras) and HD maps are integrated together to yield a more robust and accurate result. 

A Kalman filter is derived based on the assumptions that the probability distribution is Gaussian and the dynamics are linear. However, practical models generally involve nonlinear dynamics and non-Gaussian distributions, and for such applications, Kalman filters are less accurate than optimization based methods.

 There have been ongoing efforts to improve the performance of Kalman filters. The unscented Kalman filter (UKF)~\cite{ukf2000} is used to reduce the local linearization error of the extended Kalman filer (EKF); inverse depth parametrization~\cite{inversedepth2006} is used to represent the three-dimensional spatial location of feature points recognized by camera, as a Gaussian distribution on Euclidean space about a feature point's position does not model its depth uncertainty well; the multi-state constraint Kalman filter (MSCKF)~\cite{msckf2007} method only maintains a sliding window of historic camera poses, and removes feature points from the state space. Multiple views of the same feature point gives a constraint on the state variable, and serves as measurement. Such an approach can be viewed as an attempt to further relax the assumption that the uncertainty in the feature point positions are described by Gaussian distributions and thereby improve the algorithm's performance.

There is also an effort to improve the consistency of the Kalman filter by taking into account the issue of observability. The EKF does not respect the observability properties of the underlying continuous control system, which is believed to be a major source of inconsistency in the state estimate. It is hoped that by correctly incorporating the observability property into the construction of the algorithm, the consistency and hence accuracy of the method will be improved. An example of this phenomena is given in~\cite{counterexample2001}, where a stationary robot measured the relative position of new landmarks multiple times, and the covariance of robot's pose estimate becomes smaller over time. This is contrary to our intuition, as the the robot is static and therefore new landmarks do not provide new information about the robot's pose, and re-observation of these landmarks should not affect robot's own pose estimation. So, the pose of the robot is unobservable given such measurements, but the algorithm mistakenly treats it as an observable variable.

This phenomena has been studied in depth using observability analysis~\cite{huang2008analysis}. It turns out three dimensions of the state space, i.e., the robot's position and its orientation, are unobservable to odometry and camera measurements, while for the linearized discrete model, only two dimensions (the robot's position) remain unobservable. This was addressed by considering two modifications, one is the special transition matrix $\Phi(\Check{x}_{k+1}, \hat{x}_k)$ that depends on $\hat{x}_k$ (modified state value at $k$-th step) and $\Check{x}_{k+1}$ (predicted state value at $(k+1)$-th step). If we denote the unobservable subspace at $\hat{x}_k$ by $\mathcal{N}_k$, then after the measurement update, $\Phi(\Check{x}_{k+1}, \hat{x}_k)\mathcal{N}_k$ is still perpendicular to $dh(\Check{x}_{k+1})$. More importantly, this transition matrix satisfies $\Phi(x_{k+2}, x_{k+1})\cdot \Phi(x_{k+1}, x_k) = \Phi(x_{k+2},x_k)$, thus along the true trajectory, we have
\begin{equation}\label{chainproperty}
    \Phi(x_{k+s}, x_{k+s-1})\dots \Phi(x_{k+2}, x_{k+1})\cdot \Phi(x_{k+1}, x_k)\mathcal{N}_k = \Phi(x_{k+s}, x_k)\mathcal{N}_k  \perp dh(x_{k+s}),
\end{equation}
which means the unobservable subspace remains unobservable during the propagation, which is a discrete version of the invariance property satisfied by the continuous control system.  However, since the Kalman filter is implemented by alternating propagation and measurement update, the Jacobians are calculated at both $\hat{x}_k$ and $\check{x}_k$, and as a consequence, property \eqref{chainproperty} no longer holds; based on this observation, the other modification is the First Estimate Jacobian (FEJ) framework, where Jacobians are evaluated on the first ever available estimate of all the state variables so to make property \eqref{chainproperty} hold again. Later, a revised technique termed observability constrained (OC)-EKF~\cite{huang2010observability} was proposed to not only guarantee the desired observability property but also to minimize the expected errors of the linearization points. These techniques were combined with the MSCKF method to obtain MSCKF2.0~\cite{msckf2}, which is an efficient algorithm that claims to have performance no worse than optimization-based methods.  
 
In this paper, we try to approach the observability problem from a different perspective. Notice that the Kalman filter is composed of alternating propagation and measurement update steps. While the propagation step is simply a discretization of underlying stochastic differential equation, the measurement update step is essentially Bayesian inference. Such a measurement update step does not exist in the continuous control system when we perform observability analysis, so the propagation of system is ``interrupted" at each step, so forcing this discrete system to mimic the observability property of the continuous system is somewhat artificial. Since the measurement update step is essentially Bayesian inference, we would try to solve the observability problem at this step from the perspective of probability theory. The main idea is that by viewing the unobservability of the system as an invariance with respect to a group action, the state space can be decomposed into unobservable and observable parts explicitly, $x = (x_N, x_O)$. As in the case of symmetry reduction in geometric mechanics, we can obtain reduced equations that only involve the reduced (or observable) variables, and the propagation and measurement equations will be expressed only in terms of the observable variables $x_O$. We claim that in order to keep $x_N$ unobservable during the algorithm, i.e., avoid introducing spurious information from measurements to $x_N$, Bayesian inference should only be employed on the $x_O$ part, which yields a reduced measurement model. Mathematically, it means, for a probability distribution $p(x_N, x_O)$, which can be factored into
$$p(x_N, x_O) = p(x_O)\cdot p(x_N\mid x_O),$$
the conditional distribution $p(x_N\mid x_O)$ shall remain unchanged during the measurement update, and only $p(x_O)$ can be improved by the reduced measurement. So after the measurement update
\begin{equation}\label{reduced_update}
  \overline{p}(x_N, x_O) = \overline{p}(x_O)\cdot p(x_N\mid x_O),  
\end{equation}
where $\overline{p}(x_O)$ is the improved distribution for the observable variables. Even though measurement is not applied directly to $x_N$, improvement in the estimation of $x_O$ helps to improve the estimation of $x_N$ implicitly via the conditional probability $p(x_N\mid x_O)$, which is a consequence of the correlation between $x_N$ and $x_O$. A geometric picture may help to understand this, the unobservable subspaces are tangent spaces to a collection of submanifolds in state space, and $x_N$ denotes the degrees of freedom on these submanifold. Each point on a given submanifold will yield the same measurement results. The $x_O$ indicates which submanifold the state is in, and the different submanifolds are distinguishable using the measurement data. Suppose we have a probability distribution on state space, if we are more certain about $x_O$, i.e., which submanifold we are in, it also improves our estimate about $x_N$ as long as they are not independent. By using this framework, the inconsistency counterexample in \cite{counterexample2001} is naturally avoided. In that case, the unobservable variables are $x_N = (p, \theta)$, the position and orientation of robot, and the observable variables are $x_O = (z_1, z_2\dots, z_N)$, where $z_i$ is the relative position of $i$-th feature point. When new features are observed, new $z_i$ will be appended to the state variables, and they are independent from the old state variables. As long as the robot remains static, i.e., no propagation of the uncertainty is performed, no matter how often we re-observe these feature points and improve their relative position estimation, it will not affect the pose of the robot itself.

In Section \ref{background}, we review the planar robot model from~\cite{huang2010observability}, together with its observability property. We stress the use of group actions to characterize unobservable subspaces.  In Section \ref{geometric_reduction}, we discuss how to obtain a reduced control system by considering the invariant group action, and how to decompose the planar robot system into the $x_N$ and $x_O$ parts. In Section \ref{bayesian}, we choose an appropriate Riemannian structure for the new coordinate system, and implement reduced Bayesian inference \eqref{reduced_update} for the Gaussian distribution. Numerical experiments are given in Section \ref{numerical}, where straight line, circle and general trajectories are tested, and the reduced EKF method demonstrates superior performance when compared to the classical EKF and FE-EKF method.

\section{Background}\label{background}
A planar robot is considered, equipped with odometry and cameras, navigating in an unknown environment. The odometry measures local velocity and angular velocity of the robot, while cameras detect feature points and measure their relative position to the robot. To make the notation simple and the numerical simulation easier to implement, we assume that the robot is equipped with binocular cameras, thus at each instant, measurements from two cameras will be combined to give an estimate of the relative position of feature points, and this serves as our measurement. Fix a global frame, and assign an intrinsic frame $\{e_1, e_2\}$ to the robot. The robot pose is denoted by $(p, \theta)$, where $p$ is the global position of the robot, and $\theta$ represents the relative angle from the global frame to the intrinsic frame. Thus, the intrinsic frame of the robot is given by
$$\{e_1, e_2\} = R(\theta) = \begin{pmatrix}
\cos\theta & -\sin\theta \\
\sin\theta & \cos\theta
\end{pmatrix}.$$
The global position of feature points are denoted by $(p^f_1, p^f_2,\dots p^f_N)$, and together, they form the classical state space of this problem $$X = (p, \theta, p^f_1, p^f_2,\dots p^f_N).$$ The total dimension of $X$ is $3+2N$, and this state space has the structure $SE(2)\times \mathbb{R}^{2N}$, where $SE(2)$ is the two-dimensional Euclidean group. Denote the local velocity and angular velocity measured using odometry by $v$ and $w$, and the relative position of the $i$-th feature point by $z_i$. Then, the underlying continuous control system is given by
\begin{equation}\label{propagate}
\left\{
\begin{aligned}
\dot{p} &= R(\theta)v, \\
\dot{\theta} &= w, \\
\dot{p}^f_i &= 0;
\end{aligned} \right.
\end{equation}
\begin{equation}\label{measurement}
    z_i = R(\theta)^{T}(p^f_i - p),
\end{equation}
where \eqref{propagate} is the propagation equation, and \eqref{measurement} is the measurement equation. 

As we have already mentioned, there are three unobservable dimensions in this model, the robot position $p$ and its orientation $\theta$. There are two main ways to detect unobservable dimensions in a given control system. One is the observability rank criterion \cite{hermann1977nonlinear}. Given a general affine control system
\begin{equation}\label{affine_system}
\left\{
\begin{aligned}
x' &= f(x) + \sum\limits_{i=1}^ng_i(x)u_i, \\
z &= h(x),
\end{aligned} \right.
\end{equation}
regard $f(x)$, $g_1(x)$, $g_2(x)$, $\dots$, $g_n(x)$ as vector fields on the state space $X$, then we construct the minimal dual distribution that is invariant under $f(x)$, $g_1(x)$, $g_2(x)$, $\dots$, $g_n(x)$, which is given by
$$\Omega = \langle f, g_1, g_2\dots, g_n\mid dh\rangle.$$
This step involves repeated Lie derivative calculations, and is analogous to computing the observability matrix for constant linear systems. Its orthogonal distribution $\Omega^{\perp}$, which is analogous to the null space of the observability matrix, defines the unobservable subspace for this system. There exists a corresponding local coordinate transformation $\phi(x) = (x_1, x_2)$ to the distribution $\Omega^{\perp}$,  such that the system \eqref{affine_system} can be decomposed as follows
\begin{equation}\label{decomposed_system}
\left\{
\begin{aligned}
x_1' &= \overline{f}_1(x_1, x_2) + \sum\limits_{i=1}^n\overline{g}_{1i}(x_1,x_2)u_i,\\
x_2' &= \overline{f}_2(x_2) + \sum\limits_{i=1}^n\overline{g}_{2i}(x_2)u_i, \\
z &= \overline{h}(x_2). 
\end{aligned} \right.
\end{equation}
It is clear that the measurement only depends on the $x_2$ part, and $x_2$ evolves by itself, which is unaffected by $x_1$. For any solution $(x(t), u(t), z(t))$ that satisfies the system \eqref{affine_system}, if we move the initial state $x(t_0)$ along the unobservable submanifold that passes through it, i.e., by only changing the $x_1$ part, the new trajectory corresponds to the same control $u(t)$ and measurement $z(t)$, and thus $x_1$ is unobservable in this system. Any physical property of this system that is observable has to be independent of $x_1$, in other words, has to be constant along each unobservable submanifold. It has been verified in~\cite{huang2010observability} that for the system \eqref{propagate}, \eqref{measurement},
\begin{equation}\label{nullspace}
    \Omega^{\perp} = \mathcal{N} = \text{span}\begin{pmatrix}
I_2 & Jp\\
\textbf{0} & 1\\
I_2 & Jp^f_1\\
I_2 & Jp^f_2 \\
\vdots & \vdots \\
I_2 & Jp^f_N
\end{pmatrix},
\end{equation}
where 
\begin{equation*}
    J = \begin{pmatrix}
0 & -1\\
1 & 0
\end{pmatrix}.
\end{equation*}
The first two columns indicate that the robot position $p$ is unobservable, and the third column indicates that the robot orientation is also unobservable. 

The disadvantage of the observability rank criterion technique is it involves massive amounts of Lie derivative calculations. An easier way to detect unobservable dimensions is by invariant group actions. In classical mechanics, for Hamiltonian systems, we have the famous Noether's theorem, which essentially relates symmetries to conserved quantities. For example, if the system is invariant under translation, then the total linear momentum is preserved by the system; if system is invariant under rotations, then the total angular momentum is also be preserved. A similar idea can be applied to observability analysis. If our control system is invariant under some group action, i.e., corresponds to the same control $u(t)$ and measurement $z(t)$, then we have discovered an unobservable dimension of this system. For the planar robot model, we can verify that under the following translation and rotation,
\begin{equation}\label{translation_rotation}
\left\{
\begin{aligned}
p &\mapsto p + \Delta p, \\
\theta &\mapsto \theta, \\
p^f_i &\mapsto p^f_i + \Delta p;
\end{aligned} \right.\qquad
\left\{
\begin{aligned}
p &\mapsto R(\Delta\theta)p, \\
\theta &\mapsto \theta + \Delta\theta, \\
p^f_i &\mapsto R(\Delta\theta)p^f_i,
\end{aligned} \right.
\end{equation}
\eqref{propagate} \eqref{measurement} remains invariant. For each state $(p, \theta, p^f_1, p^f_2,\dots p^f_N)$, we can calculate the tangent subspace induced by the above group action. In particular, the translation action gives the first two columns of \eqref{nullspace}, and the rotation action gives the third column of \eqref{nullspace}. It is sometimes natural to find these invariant group actions by physical intuition. Odometry measures local velocity and local angular velocity, while the camera measures relative position of feature points, as such, the information that they provide is ``local" and cannot give constraints on the global pose. Global position can be measured by using GPS or if a map provides the absolute position of feature points. In the case of VIO, the IMU provides measurements of the local acceleration and angular velocity. Since the IMU measurement can detect the gravity direction, gravity breaks the full rotational symmetry, so only rotations around the gravity direction keep the system invariant, which indicates that the yaw angle is unobservable in VIO applications. The main disadvantage of the invariant group action technique as a method of characterizing unobservable dimensions is that we cannot determine if we have found all the unobservable dimensions.

We now describe the special transition matrix for \eqref{propagate} that satisfies the chain rule property \eqref{chainproperty}. Let $v = v_m + n_v$, $w = w_m + n_w$, where $v_m$,$w_m$ are odometry measurements, and $n_v$, $n_w$ are measurement noise. A simple forward Euler integration of \eqref{propagate} gives
\begin{equation}\label{propagate_noise}
\left\{
    \begin{aligned}
    p_{k+1} &= p_k + \Delta t\cdot R(\theta_k)(v_m^{k+1} + n_v), \\
    \theta_{k+1} &= \theta_k + \Delta t\cdot(w_m^{k+1} + n_w), \\
    p^f_{i,k+1} &= p^f_{i,k},
    \end{aligned}\right.
\end{equation}
so, the nominal values are updated as follows
\begin{equation}\label{discrete_propagate}
\left\{
    \begin{aligned}
    p_{k+1} &= p_k + R(\theta_k)\Delta d_{k+1}, \\
    \theta_{k+1} &= \theta_k + \Delta \theta_{k+1}, \\
    p^f_{i,k+1} &= p^f_{i,k},
    \end{aligned}\right.
\end{equation}
with $\Delta d_{k+1} = v_m^{k+1}\cdot \Delta t$, $\Delta \theta_{k+1} = w_m^{k+1}\cdot \Delta t$. The transition matrix for \eqref{discrete_propagate} is
\begin{equation}\label{transition_matrix}
    \Phi = \begin{pmatrix}
    I_2 & JR(\theta_k)\Delta d_{k+1} & \\
    \textbf{0} & 1 &\\
    & & I_{2N}
    \end{pmatrix} = \begin{pmatrix}  
    I_2 & J(p_{k+1}-p_k) & \\
    \textbf{0} & 1 &\\
    & & I_{2N}    
    \end{pmatrix} = \Phi(x_{k+1}, x_k),
\end{equation}
The derivative of the measurement \eqref{measurement} is given by
\begin{equation*}
   \begin{aligned}
    dh_i(x) &= \begin{pmatrix}-R(\theta)^T, & R(\theta)^TJ^T(p^f_i-p), &\textbf{0}, &\dots  & R(\theta)^T \dots \end{pmatrix} \\
    & = R(\theta)^T\begin{pmatrix}-I_2, & J^T(p^f_i-p), &\textbf{0}, &\dots  &  I_2\dots \end{pmatrix},
   \end{aligned}
\end{equation*}
where $I_2$ appears at the $i$-th index of the feature points. We first check that $dh_i\cdot \mathcal{N} = 0$, which is theoretically guaranteed; then the nullspace $\mathcal{N}_k$ at $x_k$ after applying the transition matrix $\Phi(x_{k+1}, x_k)$ is
\begin{equation*}
    \Phi(x_{k+1}, x_k)\cdot \mathcal{N}_k = \begin{pmatrix}  
    I_2 & J(p_{k+1}-p_k) & \\
    \textbf{0} & 1 &\\
    & & I_{2N}    
    \end{pmatrix}\cdot \begin{pmatrix}
I_2 & Jp_k\\
\textbf{0} & 1\\
I_2 & Jp^f_{1,k}\\
I_2 & Jp^f_{2,k} \\
\vdots & \vdots \\
I_2 & Jp^f_{N,k}
\end{pmatrix} = \begin{pmatrix}
I_2 & Jp_{k+1}\\
\textbf{0} & 1\\
I_2 & Jp^f_{1,k}\\
I_2 & Jp^f_{2,k} \\
\vdots & \vdots \\
I_2 & Jp^f_{N,k}
\end{pmatrix},
\end{equation*}
it is easy to verify that $dh_i(x_{k+1})\cdot \Phi(x_{k+1}, x_k)\cdot \mathcal{N}_k = 0$, since $p^f_{i,k+1} = p^f_{i,k}$. Finally, we can also verify that transition matrix $\Phi(x_{k+1}, x_k)$ \eqref{transition_matrix} satisfies $\Phi(x_{k+2}, x_k) = \Phi(x_{k+2}, x_{k+1}) \cdot \Phi(x_{k+1}, x_k)$,
so the chain rule property \eqref{chainproperty} is satisfied along the exact trajectory.

\section{Geometric reduction}\label{geometric_reduction}
Given the state space $X$, and a general control system
\begin{equation}\label{general_control}
\left\{
\begin{aligned}
x' &= f(x,u), \\
z &= h(x),
\end{aligned}\right.
\end{equation}
suppose we have a left Lie group action of $G$ on $X$, i.e., $G\times X\to X$, which keeps the system \eqref{general_control} invariant, i.e., for any $(x(t), u(t), z(t))$ that satisfies \eqref{general_control}, $(g\cdot x(t), u(t), z(t))$ also satisfies \eqref{general_control} for $\forall g\in G$, i.e.,
\begin{equation}\label{g_general_control}
\left\{
\begin{aligned}
(g\cdot x)' &= f(g\cdot x,u), \\
z &= h(g\cdot x).
\end{aligned}\right.
\end{equation}
For each point $x\in X$, $\forall \xi\in \mathfrak{g}$, where $\mathfrak{g}$ is the Lie algebra of $G$, the infinitesimal generator $\xi_X(x)=\left.\frac{d}{dt}\right|_{t=0}\text{exp}(\xi(t))\cdot x$ gives one unobservable direction of $x$, and the orbit of $x$,
$$\text{Orb}(x) = \{g\cdot x\mid g\in G\},$$
gives the unobservable submanifold that passes through $x$. This invariant group action on the control system allows us to perform reduction to express the system in terms of reduced variables on the quotient space $X/G$. This kind of reduction due to the presence of a continuous symmetry group arises in many disciplines. In optimization, when objective function is invariant under a group action, the problem can be rephrased as an optimization problem on the quotient space~\cite{absil2009optimization}; also for classical mechanics on a Lie group, when Lagrangian or Hamiltonian is invariant under the group action, the mechanics can be reduced to the Lie algebra $\mathfrak{g}$ or its dual $\mathfrak{g}^*$, which is referred to as Euler--Poincar\'{e} reduction and Lie--Poisson reduction~\cite{montaldi2005geometric}, respectively. Here, suppose that the group action $G\times X\to X$ is free and proper, then we obtain a smooth quotient space $X/G$, with a quotient map\cite{lee2001introduction}
$$\pi : X\to X/G$$
that is a smooth submersion. Then, the control system~\eqref{general_control} can be reduced to the quotient space $X/G$. For $[x]\in X/G$,
\begin{equation}\label{reduced_system}
\left\{
\begin{aligned}
[x]' &= f([x],u), \\
z &= h([x]).
\end{aligned}\right.
\end{equation}
\eqref{reduced_system} is a reduced control system with reduced propagation and measurement equations. Consider \eqref{propagate}, \eqref{measurement} as a concrete example. The state space is $X = SE(2)\times \mathbb{R}^{2N}$, and  we already know that the group action $G = SE(2)$ \eqref{translation_rotation} on $X$ leaves \eqref{propagate} \eqref{measurement} invariant, so this induces a reduced control system on $X/G$. However, $X/G$ is an abstract quotient manifold, so in order to deal with it explicitly, we need a concrete coordinate representation. It turns out that the relative position of feature points $z_i$ provide a natural coordinate representation for $X/G$. With these coordinates for the reduced space, the measurement equation \eqref{measurement} reduces to 
$$z_i = \textbf{Id}(z_i);$$
and the propagation equation \eqref{propagate} reduces to
\begin{equation*}
\begin{aligned}
    \frac{d}{dt} z_i &= \frac{d}{dt} R(\theta)^T(p^f_i-p) \\
    & = R(\theta)^{T}(\dot{p}^f_i - \dot{p}) + J^{T}R(\theta)^{T}\dot{\theta}(p^f_i-p) \\
    & = R(\theta)^T(-R(\theta)v) + J^TR(\theta)^Tw(p^f_i-p) \\
    & = -v - wJz_i.
\end{aligned}
\end{equation*}
Together with the robot position $p$ and its orientation $\theta$, we obtain a new coordinate representation of the state space, $(p, \theta, z_1, z_2\dots, z_N)$, and the control system with respect to this new coordinate system is given by
\begin{equation}\label{reduced_propagate}
\left\{
\begin{aligned}
\dot{p} &= R(\theta)v, \\
\dot{\theta} &= w, \\
\dot{z_i} &= -v - wJz_i,
\end{aligned} \right.
\end{equation}
\begin{equation}\label{reduced_measurement}
    z_i = \textbf{Id}(z_i).
\end{equation}
The coordinate system given by $(p, \theta, z_1, z_2\dots, z_N)$ is global, and can be regarded as the state space being decomposed into the product of $x_N = (p, \theta)$ and $x_O = (z_1, z_2\dots, z_N)$, the former is unobservable part, and the latter is observable part, while the group action \eqref{translation_rotation} now acts trivially on the $x_O$ part. In this case, the unobservable subspace at each point can be represented as
\begin{equation}\label{new_nullspace}
    \Omega^{\perp} = \mathcal{N} =  \text{span}\begin{pmatrix}
I_2 & Jp\\
\textbf{0} & 1\\
\textbf{0} & 0\\
\textbf{0} & 0 \\
\vdots & \vdots \\
\textbf{0} & 0
\end{pmatrix}.
\end{equation}
We improved the linearity of the measurement equation~\eqref{reduced_measurement} by transforming to the relative feature position representation, and as a consequence, the originally trivial propagation equation $\dot{p}^f_i = 0$ now becomes the nontrivial $\dot{z}_i = -v - wJz_i$. 

We can construct a similar transition matrix for \eqref{reduced_propagate} as in \eqref{transition_matrix}, and apply forward Euler integration to \eqref{reduced_propagate}, which yields

\begin{equation}\label{new_propagate_noise}
\left\{
\begin{aligned}
    p_{k+1} &= p_k + \Delta t\cdot R(\theta_k) (v_m^{k+1} + n_v), \\
    \theta_{k+1} &= \theta_k + \Delta t\cdot (w_m^{k+1} + n_w), \\
    z_{i,k+1} &= z_{i,k} + \Delta t\cdot(-(v_m^{k+1}+n_v)-(w_m^{k+1}+n_w)Jz_{i,k}),
\end{aligned} \right.
\end{equation}
and the corresponding nominal value update step is given by
\begin{equation}
    \left\{
    \begin{aligned}
        p_{k+1} &= p_k + R(\theta_k)\Delta d_{k+1}, \\
        \theta_{k+1} &= \theta_k + \Delta \theta_{k+1}, \\
        z_{i,k+1} &= (I - \Delta \theta_{k+1}J)z_{i,k} - \Delta d_{k+1},
    \end{aligned}\right.
\end{equation}
and the transition matrix is
\begin{equation}\label{new_transition_matrix}
    \Phi_{k+1,k} = \begin{pmatrix}
    I_2 & J(p_{k+1}-p_k) & & & & \\
    \textbf{0} & 1 & & & & \\
    & & I-\Delta \theta_{k+1}J & & & \\
    & & & I-\Delta \theta_{k+1}J & & \\
    & & & & \ddots & \\
    & & & & & I-\Delta \theta_{k+1}J
    \end{pmatrix}.
\end{equation}
The derivative of the reduced measurement \eqref{reduced_measurement} is
\begin{equation*}
    dh_i = \begin{pmatrix} \textbf{0}_{2\times 2}, & \textbf{0}_{2\times 1}, & \textbf{0}_{2\times 2} & \dots & I_2 & \dots  \end{pmatrix},
\end{equation*}
and we can see that $dh_i\cdot \mathcal{N} = 0$, which holds simply because the lower $2N\times 3$ part of $\mathcal{N}$ \eqref{new_nullspace} is zero. This property still holds after repeated measurement updates
$$\Phi_{k+s,k+s-1}\cdot \dots \Phi_{k+2,k+1}\cdot \Phi_{k+1,k}\cdot \mathcal{N}_k$$
as the bottom left parts of the transition matrices $\Phi_{k+j,k+j-1}$ vanish.

\section{Bayesian inference}\label{bayesian}
In Section \ref{geometric_reduction}, we discussed how to obtain the reduced control system that arises as a consequence of the invariant group action, which lead to the reduced propagation and measurement equations on the reduced quotient space. This motivated the use of a new coordinate representation $(p,\theta. z_1,z_2\dots z_N)$ for the planar robot system. This decomposes the state space into the product of the unobservable part $(p,\theta)$, and the observable part $(z_1,z_2\dots z_N)$. We also constructed a transition matrix which ensures that the unobservable subspace remains perpendicular to the measurement $dh$ during propagation. But, as we noted, this propagation is ``interrupted" at each step during the Kalman update. In order to preserve the unobservability property of the system during the measurement update, we propose a reduced Bayesian inference update\eqref{reduced_update} on the observable part using the reduced measurement. It is natural to apply Bayesian inference on $p(z_1,z_2\dots z_N)$, and update $p(p,\theta,z_1,z_2\dots z_N)$ using \eqref{reduced_update}. However, before we do that, there are some basic things that we need to make clear. The first question is, what is Bayesian inference? As we all know, Bayesian inference is
$$p(x\mid y) = \frac{p(x)\cdot p(y\mid x)}{p(y)},$$
which arises naturally from the fact that
\begin{equation}\label{conditional_probability}
 p(x,y) = p(x)\cdot p(y\mid x) = p(y)\cdot p(x\mid y).   
\end{equation}
where $p(x)$, $p(y)$, $p(x,y)$ are probability density functions. The next question is, what is a probability density function? As we all know, a probability distribution is a probability measure on the state space, and the probability density function is a concrete way to represent it,
\begin{equation}\label{density_definition}
    P(x\in A) = \int_A p(x)dx.
\end{equation}
However, in \eqref{density_definition} we need to specify a measure $dx$ to perform the integration, and the representation of the probability measure as a probability density function depends on that choice of $dx$. Also, given a probability distribution on state space, when we try to find the most probable point, we find this problem also requires the introduction of a measure on the state space. As such, we need to specify a measure on the state space in order to define the probability density function or to find the most probable point. For a more detailed discussion of such issues, see~\cite{jermyn2005invariant}. After specifying the measure $dx$, the density function associated with a probability measure $\mu$ is just the Radon--Nikodym derivative of $\mu$ with respect to $dx$. A natural way to specify a measure on a smooth manifold is by specifying a Riemannian structure on it. Each Riemannian structure induces a Riemannian volume form, which in turn induces a measure. When we write down the common Gaussian density function, we are actually assuming the standard Riemannian structure on $\mathbb{R}^n$ implicitly. This structure is quite natural, as it is homogeneous and isotropic and does not introduce prior information on the space. More precisely, this structure is invariant under Euclidean transformations.  

With the observation that each density function is defined with respect to an underlying base measure, we look at \eqref{conditional_probability} again, and find that there is a product measure implicitly defined on $(x,y)$, and that \eqref{conditional_probability} is essentially an application of Fubini's theorem.

Now consider the new coordinate system $(p,\theta, z_1, z_2\dots z_N)$ we get for the planar robot system, which is diffeomorphic to the original coordinate system by the following transformations
\begin{equation}\label{change_of_coordinate}
(p,\theta, p^f_1, p^f_2\dots p^f_N)\leftrightarrow (p,\theta, z_1, z_2\dots z_N): 
\left\{
\begin{aligned}
p &= p,\\
\theta &= \theta,\\
p^f_i &= p + R(\theta)z_i;
\end{aligned} \right.\qquad
\left\{
\begin{aligned}
p &= p,\\
\theta &= \theta,\\
z_i &= R(\theta)^T(p^f_i-p).
\end{aligned} \right.
\end{equation}
Given the decomposition $p(x_N,x_O) = p(x_O)\cdot p(x_N\mid x_O)$, where $x_N = (p,\theta)$, $x_O = (z_1, z_2,\dots z_N)$, which Riemannian structure should we assign to the state space to obtain a base measure? For the old coordinate system $(p,\theta, p^f_1, p^f_2\dots p^f_N)$, it is natural to use the standard Riemannian structure $dp\otimes dp + d\theta\otimes d\theta + \sum_{i=1}^Ndp^f_i\otimes dp^f_i$, as it is invariant under the group action \eqref{translation_rotation}. Moreover, if a Lie group $G$ acts freely and properly on a Riemannian manifold $M$, and the action is isometric for $\forall g\in G$, then there is a natural Riemannian structure on $M/G$ induced from $M$. For our quotient space $(z_1,z_2\dots z_N)$, it can be verified that the induced Riemannian structure is simply the standard Riemannian structure of $\mathbb{R}^{2N}$. However, the Riemannian structure of $(p,\theta,z_1,z_2\dots z_N)$ induced by the diffeomorphism \eqref{change_of_coordinate} is not simply the product of the standard Riemannian structure of $(p,\theta)$ with the standard Riemannian structure of $(z_1,z_2\dots z_N)$. This can be verified by checking the Jacobian of \eqref{change_of_coordinate}, which is not a unitary matrix. We now have to decide which Riemannian structure to use on $(p,\theta,z_1,z_2\dots z_N)$, either choose the non-product structure induced from \eqref{change_of_coordinate}, or choose the product structure. We choose the product structure here, one reason is that it is simpler to apply the reduced Bayesian inference step~\eqref{reduced_update}, another reason is this product structure is also invariant under the group action \eqref{translation_rotation}, which ensures that $(z_1,z_2\dots z_N)$ remains fixed. 

Finally, we are in a position to discuss how to implement \eqref{reduced_update} for the Gaussian distribution with respect to the standard Riemannian structure. Consider $p(x_N,x_O)$ with a Gaussian distribution given by
$$p(x_N,x_O) = \mathcal{N}\left(\begin{pmatrix}\mu_N \\ \mu_O\end{pmatrix}, \begin{pmatrix}\Sigma_{NN} &\Sigma_{NO}\\ \Sigma_{ON} & \Sigma_{OO}\end{pmatrix}\right).$$
Then, by the property of conditional Gaussian distributions,
\begin{equation*}
    \begin{aligned}
    p(x_N,x_O) &= p(x_O)\cdot p(x_N\mid x_O) \\
    & = \mathcal{N}(x_O\mid \mu_O, \Sigma_{OO})\cdot \mathcal{N}(x_N\mid \mu_{N\mid O}, \Sigma_{N\mid O}), \\
    &= \mathcal{N}(x_O\mid \mu_O, \Sigma_{OO})\cdot \mathcal{N}(x_N\mid \mu_N+\Sigma_{NO}\Sigma^{-1}_{OO}(x_O-\mu_O), \Sigma_{NN} - \Sigma_{NO}\Sigma^{-1}_{OO}\Sigma_{ON}).
    \end{aligned}
\end{equation*}
The next step is to obtain an improved $\overline{p}(x_O) = \mathcal{N}(\overline{\mu}_O, \overline{\Sigma}_{OO})$ by taking a reduced measurement and computing the improved joint distribution $\overline{p}(x_N,x_O)$. In order to compute the joint distribution from the product of two Gaussian distributions efficiently, we choose the precision matrix, information vector representation of the Gaussian distribution,
$$\mathcal{N}(x_O\mid \mu_O,\Sigma_{OO}) = \mathcal{N}(x_O\mid \xi_O, \Omega_{OO}),$$
where $\Omega_{OO} = \Sigma^{-1}_{OO}$, $\xi_O = \Sigma^{-1}_{OO}\cdot\mu_O$, and $\Omega_{N\mid O} = \Sigma^{-1}_{N\mid O}$. Suppose that the reduced discrete measurement equation is given by
$$z = C\cdot x_O + \delta,$$
with measurement noise $\delta\sim \mathcal{N}(0, Q)$, where $Q$ is the noise covariance, then by the information filter \cite{thrun2005probabilistic},
\begin{equation*}
    \left\{
    \begin{aligned}
    \overline{\Omega}_{OO} &= \Omega_{OO} + C^TQ^{-1}C, \\
    \overline{\xi}_O &= \xi_O + C^TQ^{-1}z.
    \end{aligned}\right.
\end{equation*}
Let $\mu_N+\Sigma_{NO}\Sigma^{-1}_{OO}(x_O-\mu_O) = Ax_O +b$ for notational simplicity, where $A = \Sigma_{NO}\Sigma^{-1}_{OO}$, $b = \mu_N - \Sigma_{NO}\Sigma^{-1}_{OO}\mu_O$. We calculate the density function of $\overline{p}(x_N,x_O)$,
\begin{equation*}
    \begin{aligned}
    \text{log}\overline{p}(x_N,x_O) &= \text{log}\overline{p}(x_O) + \text{log}p(x_N\mid x_O) \\
    &= \text{log}\mathcal{N}(x_O\mid \overline{\xi}_O, \overline{\Omega}_{OO}) + \text{log}\mathcal{N}(x_N\mid A\cdot x_O+b, \Omega_{N\mid O})\\
    & \simeq -\frac{1}{2}x_O^T\overline{\Omega}_{OO}x_O + \overline{\xi}_O^Tx_O - \frac{1}{2}(x_N - Ax_O-b)^T\Omega_{N\mid O}(x_N - Ax_O-b)\\
    & \simeq - \frac{1}{2}x_O^T\overline{\Omega}_{OO}x_O + \overline{\xi}_O^Tx_O - \frac{1}{2}x_N^T\Omega_{N\mid O}x_N - \frac{1}{2}x_OA^T\Omega_{N\mid O}Ax_O\\
    &\qquad\quad+ x_O^TA^T\Omega_{N\mid O}x_N + x^T_N\Omega_{N\mid O}b - x_O^TA^T\Omega_{N\mid O}b\\
    & = -\frac{1}{2}(x_N,x_O)^T\begin{pmatrix}\Omega_{N\mid O} & -\Omega_{N\mid O}A\\ -A^T\Omega_{N\mid O} & \overline{\Omega}_{OO} + A^T\Omega_{N\mid O}A   \end{pmatrix}\begin{pmatrix} x_N \\ x_O \end{pmatrix}+ (x_N,x_O)^T\begin{pmatrix} \Omega_{N\mid O}b \\ \overline{\xi}_O - A^T\Omega_{N\mid O}b \end{pmatrix}.
    \end{aligned}
\end{equation*}
Thus, the improved joint probability distribution $\overline{p}(x_N,x_O)$ has the following improved precision matrix and information vector,
\begin{equation}\label{improved_precision_information}
    \overline{\Omega} = \begin{pmatrix}\Omega_{N\mid O} & -\Omega_{N\mid O}A\\ -A^T\Omega_{N\mid O} & \overline{\Omega}_{OO} + A^T\Omega_{N\mid O}A   \end{pmatrix}, \qquad \overline{\xi} = \begin{pmatrix} \Omega_{N\mid O}b \\ \overline{\xi}_O - A^T\Omega_{N\mid O}b \end{pmatrix}.
\end{equation}
We can see from \eqref{improved_precision_information} that if originally $\Sigma_{NO} = 0$, then $\Sigma_{N\mid O} = \Sigma_{NN}$ and $A = \Sigma_{NO}\Sigma^{-1}_{OO} = 0$, thus $\overline{\Omega}$ is also block diagonal, and the improvement of $\overline{\Omega}_{OO}$ and $\overline{\xi}_O$ does not affect the distribution of $x_N$.

\section{Numerical experiment}\label{numerical}
We performed numerical experiments for the planar robot system \eqref{propagate}, \eqref{measurement}, where the robot is assumed to be round with diameter $0.5\,\rm{m}$, and binocular cameras are equipped on the left and right sides of the robot to detect feature points, each with a field angle of $120^{\circ}$. The odometry measures local velocity $v$ (units: m/s) and angular velocity $w$ (units: rad) with Gaussian noise $\mathcal{N}(\textbf{0}, Q_v)$, $\mathcal{N}(0,Q_w)$, respectively. We set $Q_v = \begin{pmatrix} 0.01 & 0\\ 0 & 0.01\end{pmatrix}$ and $Q_w = 0.01$ during the simulation. The camera is assumed to measure the relative direction of the feature point with Gaussian noise $\mathcal{N}(0,Q_z)$. As we mentioned in Section \ref{background}, at each step, the measurement of feature points by the binocular cameras shall be combined together to give an estimate of the relative position of the feature points, and the uncertainty in the relative position is assumed to be described by an approximate Gaussian distribution. The odometry and camera readings are updated at the same frequency.

We tested three different methods, the first is the classical EKF on the state space $(p, \theta, p^f_1, p^f_2\dots p^f_N)$, with the propagation equation \eqref{propagate_noise}; the second is the First Estimate EKF (FE-EKF) where the First Estimate technique in \cite{huang2008analysis} is used; the third is our reduced EKF on the state space $(p,\theta, z_1,z_2\dots z_N)$, with propagation equation \eqref{new_propagate_noise} and reduced measurement update equation \eqref{improved_precision_information}. We observed superior performance of the reduced EKF compared to classical EKF and FE-EKF, especially at estimating the robot's orientation $\theta$. Furthermore, the reduced EKF is less sensitive to measurement noise.

We considered three different trajectories: straight line, circle, and a general trajectory. Feature points are generated along the  trajectory randomly with a given density. We assume that there are no errors introduced during feature detection and matching. We ignore feature points that are far away (distance $>5\,\rm{m}$) from the robot, since triangulation for such points by binocular cameras are unstable. 
\subsection{Straight line}
We tested on a straight line trajectory that is $60$ m long, and the robot traveled along it with constant speed $1\,\rm{m/s}$. We also compared low and high densities of feature points.
\begin{figure}[H]
\centering
	\begin{subfigure}[b]{0.45\textwidth}
		\includegraphics[scale=0.4]{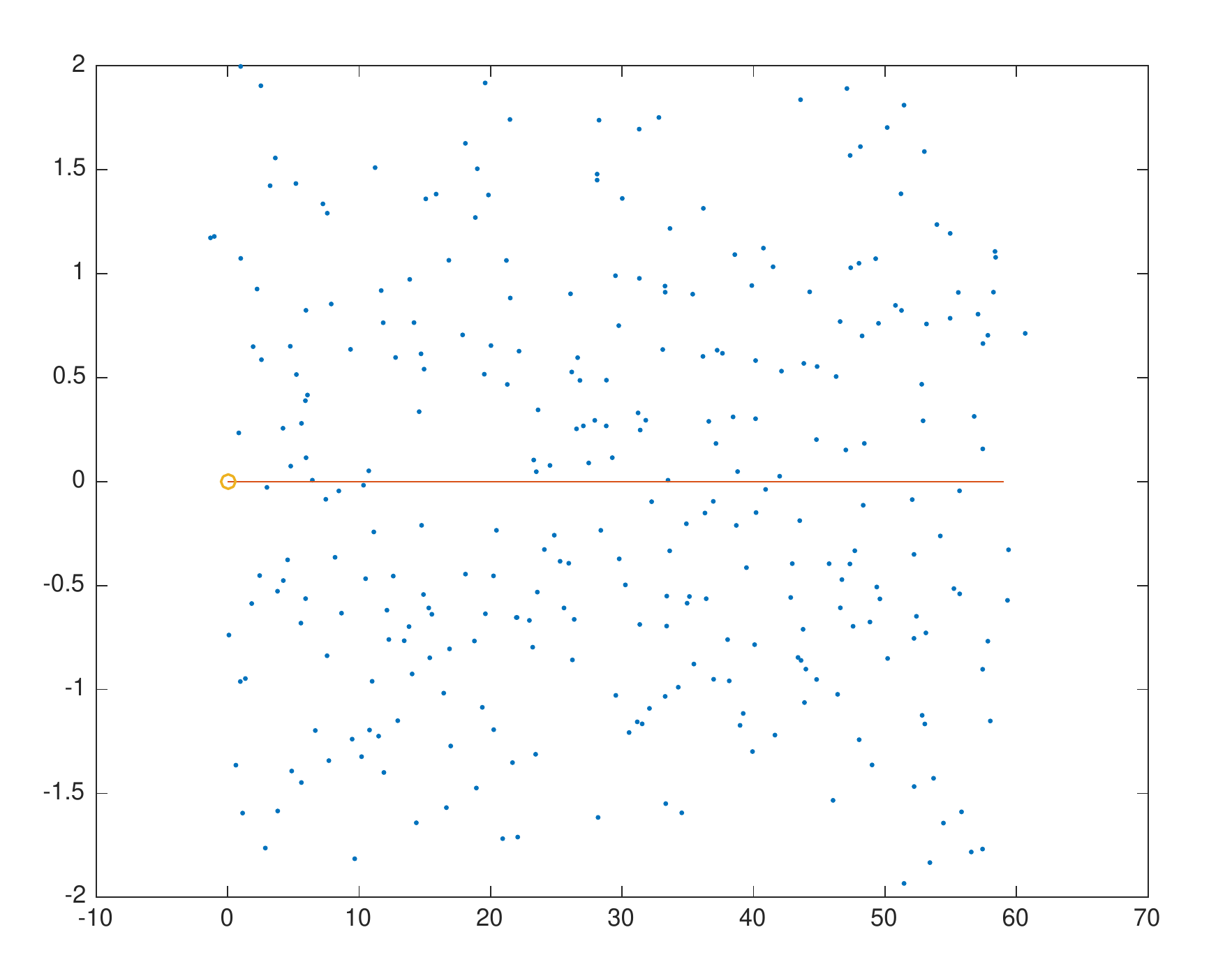}
		\caption{Straight line, low density feature distribution}
	\end{subfigure}
	\begin{subfigure}[b]{0.45\textwidth}
		\includegraphics[scale=0.4]{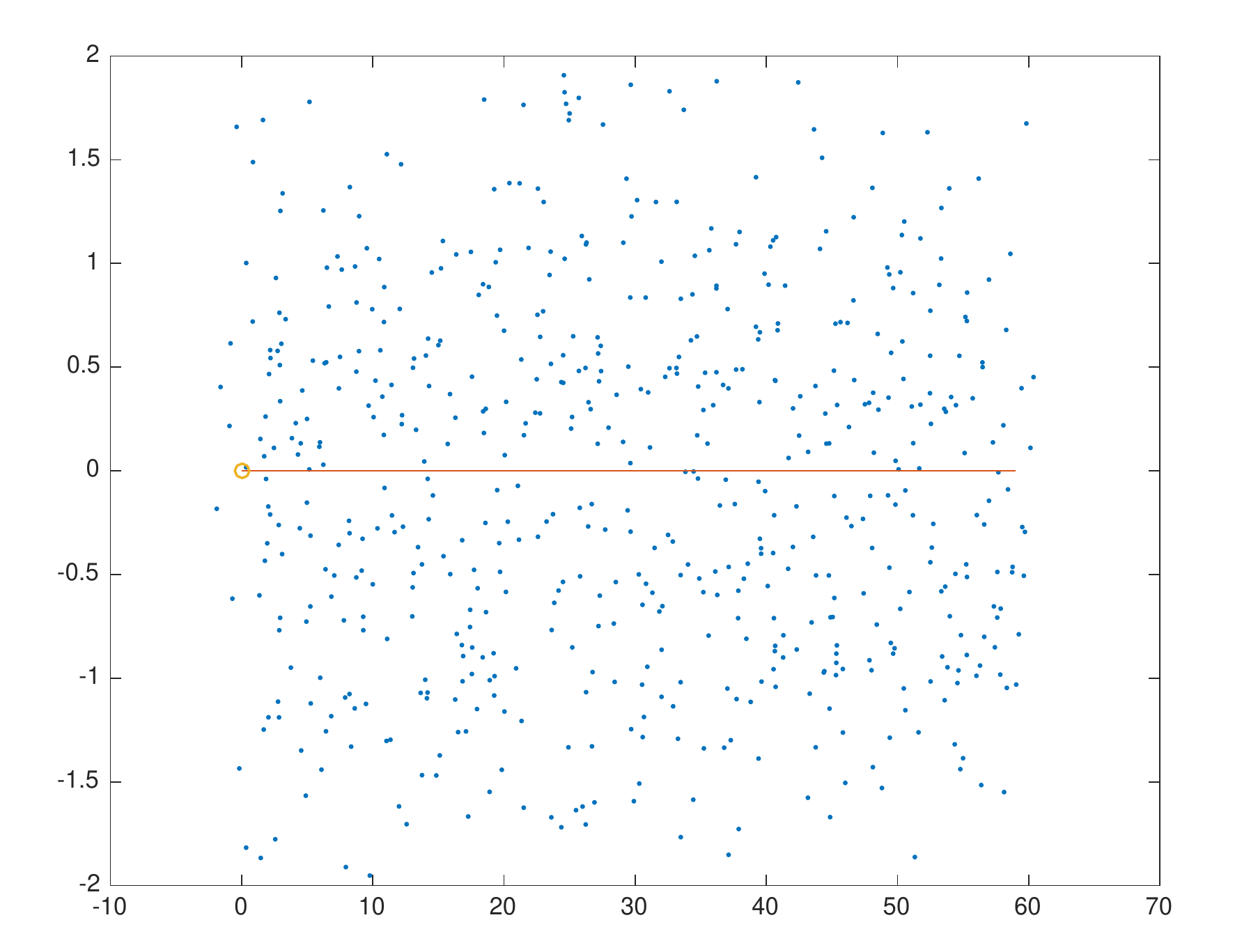}
		\caption{Straight line, high density feature distribution}
\end{subfigure}
\caption{\label{straight_traj}Straight line trajectory}
\end{figure}
In Figure~\ref{straight_traj}, the red line is the trajectory, the blue points are features, and yellow circle is the robot. For the low density feature distribution, the average feature points recognized at each step is approximately 20, and for the high density feature distribution, the average recognized feature points at each step is approximately 40. Besides the feature distribution, we also considered update frequencies of $10\,\rm{Hz}$ and $20\,\rm{Hz}$. Moreover, for fixed feature density and update frequency, we adjusted $Q_z$ to observe the sensitivity of the different methods to changes in the measurement noise. Since the feature distribution is generated randomly for fixed feature density, and odometry and camera readings are also generated randomly from the given parameters, the results for a given method will vary between realizations. The three methods (EKF, FE-EKF, reduced EKF) are applied to the same data set of feature distribution and sensor readings, and we generate 20 realizations per combination of update frequency, feature distribution, and measurement noise. We summarize the average position error (units: meters) and the orientation error (units: rad) in Table~\ref{table1} and Table~\ref{table2}.

\begin{table}\footnotesize
\centering
\begin{tabular}{|c|c|c|c|c|c|c|c|c|c|c|c|c|}
\hline
&  \multicolumn{12}{|c|}{Straight line trajectory, average feature num = 20} \\
\hline Frequency & \multicolumn{6}{|c|}{update frequency = 10Hz} & \multicolumn{6}{|c|}{update frequency = 20Hz} \\
\hline Method & \multicolumn{2}{|c|}{EKF} & \multicolumn{2}{|c|}{FE-EKF} & \multicolumn{2}{|c|}{Reduced EKF} & \multicolumn{2}{|c|}{EKF} & \multicolumn{2}{|c|}{FE-EKF} & \multicolumn{2}{|c|}{Reduced EKF} \\
\hline Error & $\delta{p}$ & $\delta{\theta}$ & $\delta{p}$ & $\delta{\theta}$ & $\delta{p}$ & $\delta{\theta}$ & $\delta{p}$ & $\delta{\theta}$ & $\delta{p}$ & $\delta{\theta}$ & $\delta{p}$ & $\delta{\theta}$ \\
\hline $Q_z=1\times 10^{-4}$ & 1.403 & 0.0375 & 1.539 & 0.0421 & 1.436 & 0.0373& 1.365 & 0.0308 & 1.456 & 0.0345 & 1.300 & 0.0282 \\
\hline $Q_z=2\times 10^{-4}$ & 2.492 & 0.0487 & 2.525 & 0.0552 & 2.549 & 0.0506 & 2.395 & 0.0411 & 2.621 & 0.057 & 2.397 & 0.0432 \\
\hline $Q_z=4\times 10^{-4}$ & 4.332 & 0.0752 & 4.284 & 0.0758 & 4.449 & 0.0773 & 4.265 & 0.0578 & 4.283 & 0.0677 & 4.317 & 0.0644 \\
\hline
\end{tabular}
\caption{\label{table1}Average error, straight line trajectory, average feature num = 20}
\end{table}

\begin{figure}[H]
\centering
	\begin{subfigure}[b]{0.45\textwidth}
		\includegraphics[scale=0.4]{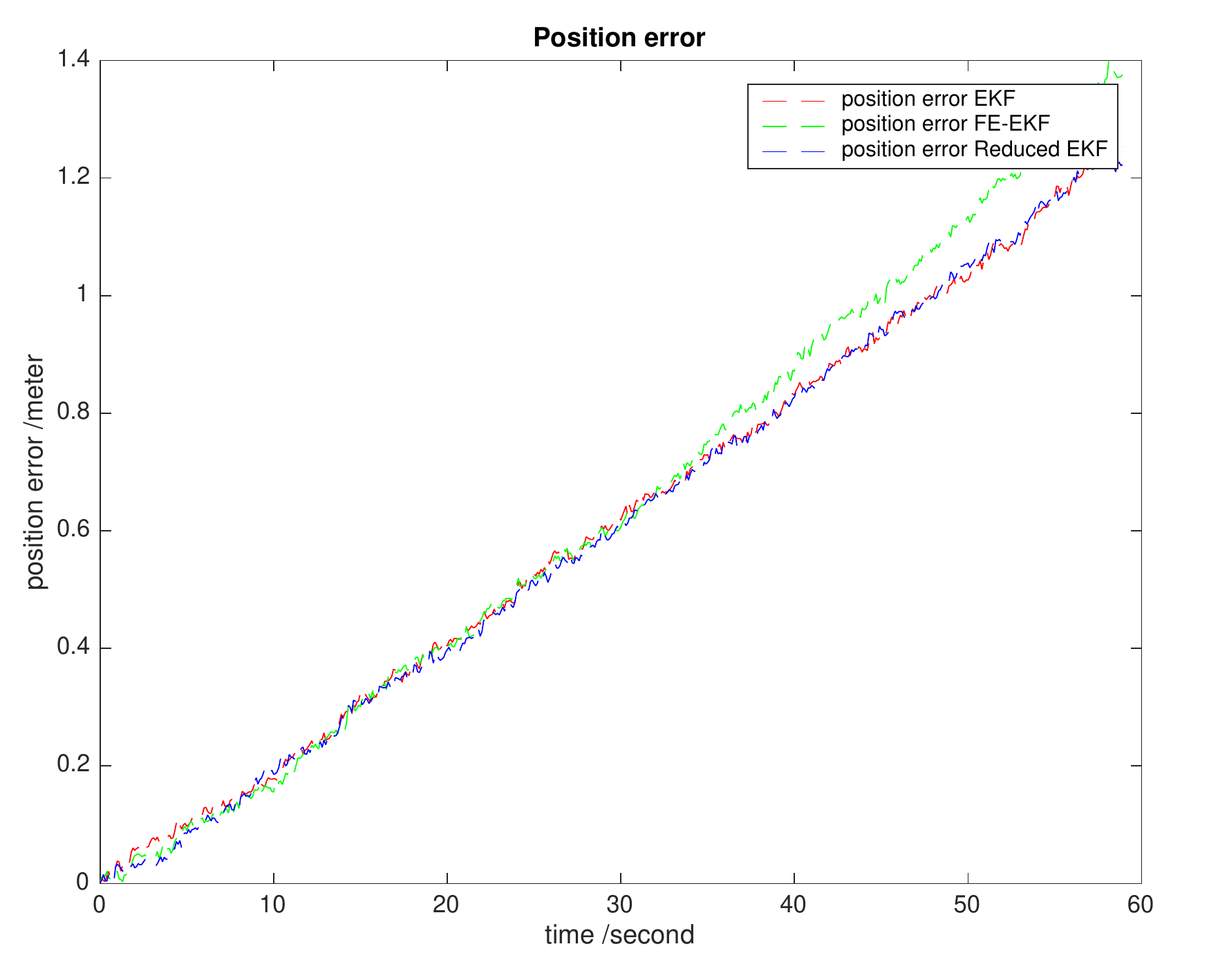}
		\caption{Position error, straight line}
	\end{subfigure}
	\begin{subfigure}[b]{0.45\textwidth}
		\includegraphics[scale=0.4]{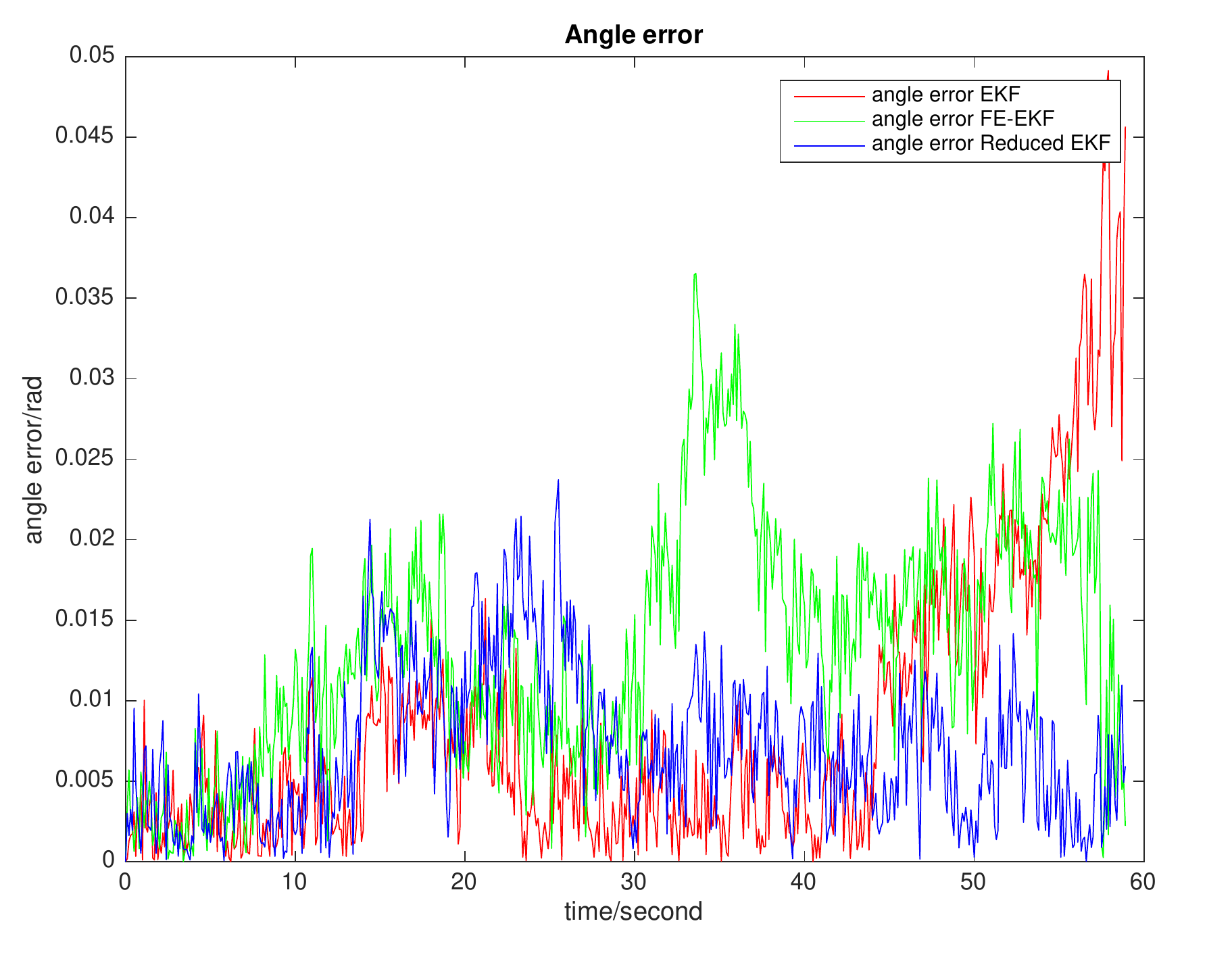}
		\caption{Angle error, straight line}
\end{subfigure}
\caption{\label{fig:error_straight}Error growth, straight line trajectory}
\end{figure}

\begin{table}\footnotesize
\centering
\begin{tabular}{|c|c|c|c|c|c|c|c|c|c|c|c|c|}
\hline
&  \multicolumn{12}{|c|}{Straight line trajectory, average feature num = 40} \\
\hline Frequency & \multicolumn{6}{|c|}{update frequency = 10Hz} & \multicolumn{6}{|c|}{update frequency = 20Hz} \\
\hline Method & \multicolumn{2}{|c|}{EKF} & \multicolumn{2}{|c|}{FE-EKF} & \multicolumn{2}{|c|}{Reduced EKF} & \multicolumn{2}{|c|}{EKF} & \multicolumn{2}{|c|}{FE-EKF} & \multicolumn{2}{|c|}{Reduced EKF} \\
\hline Error & $\delta{p}$ & $\delta{\theta}$ & $\delta{p}$ & $\delta{\theta}$ & $\delta{p}$ & $\delta{\theta}$ & $\delta{p}$ & $\delta{\theta}$ & $\delta{p}$ & $\delta{\theta}$ & $\delta{p}$ & $\delta{\theta}$ \\
\hline $Q_z=1\times 10^{-4}$ &  1.303 & 0.0269 & 1.306 & 0.027 & 1.256 & 0.0227 & 1.21 & 0.0167 & 1.283 & 0.0238 & 1.258 & 0.0208\\
\hline $Q_z=2\times 10^{-4}$ &  2.409 & 0.0385 & 2.459 & 0.0446 & 2.43 & 0.0347 & 2.331 & 0.0279 & 2.318 & 0.0313 & 2.386 & 0.0345\\
\hline $Q_z=4\times 10^{-4}$ &  4.238 & 0.047 & 4.361 & 0.0628 & 4.386 & 0.0613 & 4.23 & 0.0445 & 4.353 & 0.0558 & 4.352 & 0.0577\\
\hline
\end{tabular}
\caption{\label{table2}Average error, straight line trajectory, average feature num = 40}
\end{table}
We see that the accuracy of all three methods improve a little bit when we double the update frequency and increase the average feature number per step, which is quite natural. In addition, when we increase the uncertainty of measurement $Q_z$, the performance of all three methods degrade. There is no evident improvement in accuracy for both FE-EKF and reduced EKF methods over the classical EKF method for the straight line case. A typical error growth in the position and orientation for the straight line trajectory is given in Figure~\ref{fig:error_straight}.

\subsection{Circle}
\begin{figure}[H]
\centering
	\begin{subfigure}[b]{0.45\textwidth}
		\includegraphics[scale=0.4]{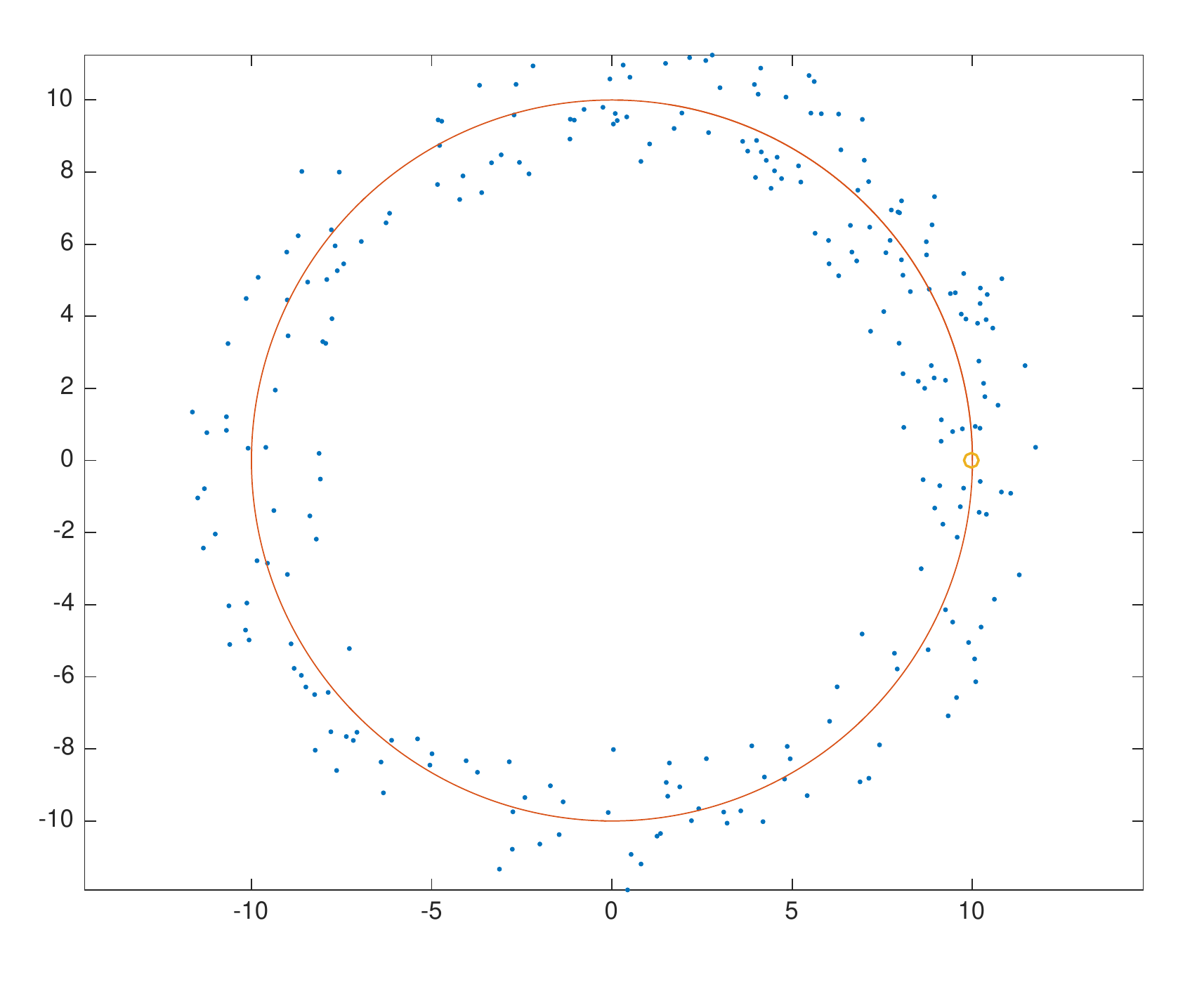}
		\caption{Circle, low density feature distribution}
	\end{subfigure}
	\begin{subfigure}[b]{0.45\textwidth}
		\includegraphics[scale=0.4]{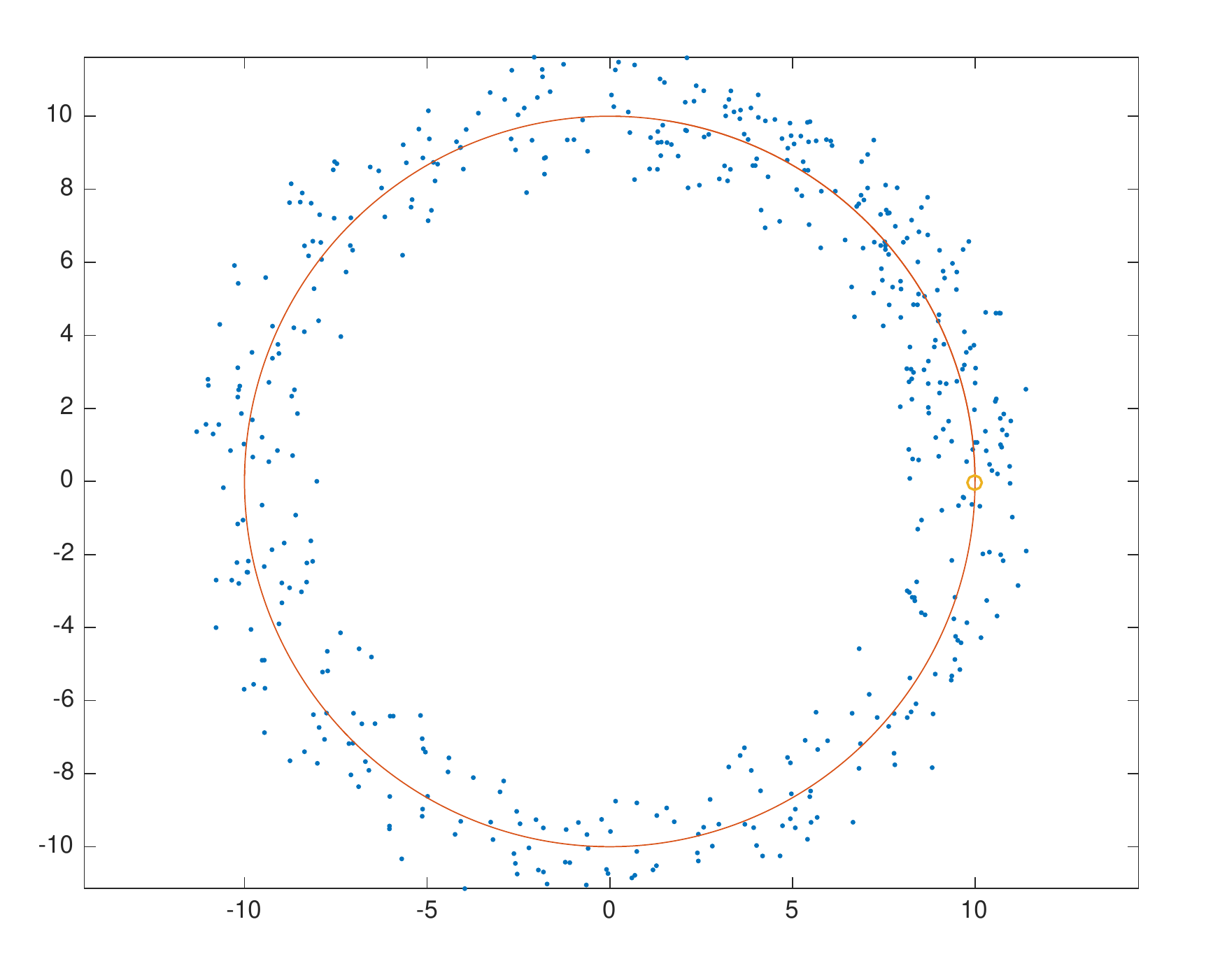}
		\caption{Circle, high density feature distribution}
\end{subfigure}
\begin{subfigure}[b]{0.45\textwidth}
   \includegraphics[scale=0.4]{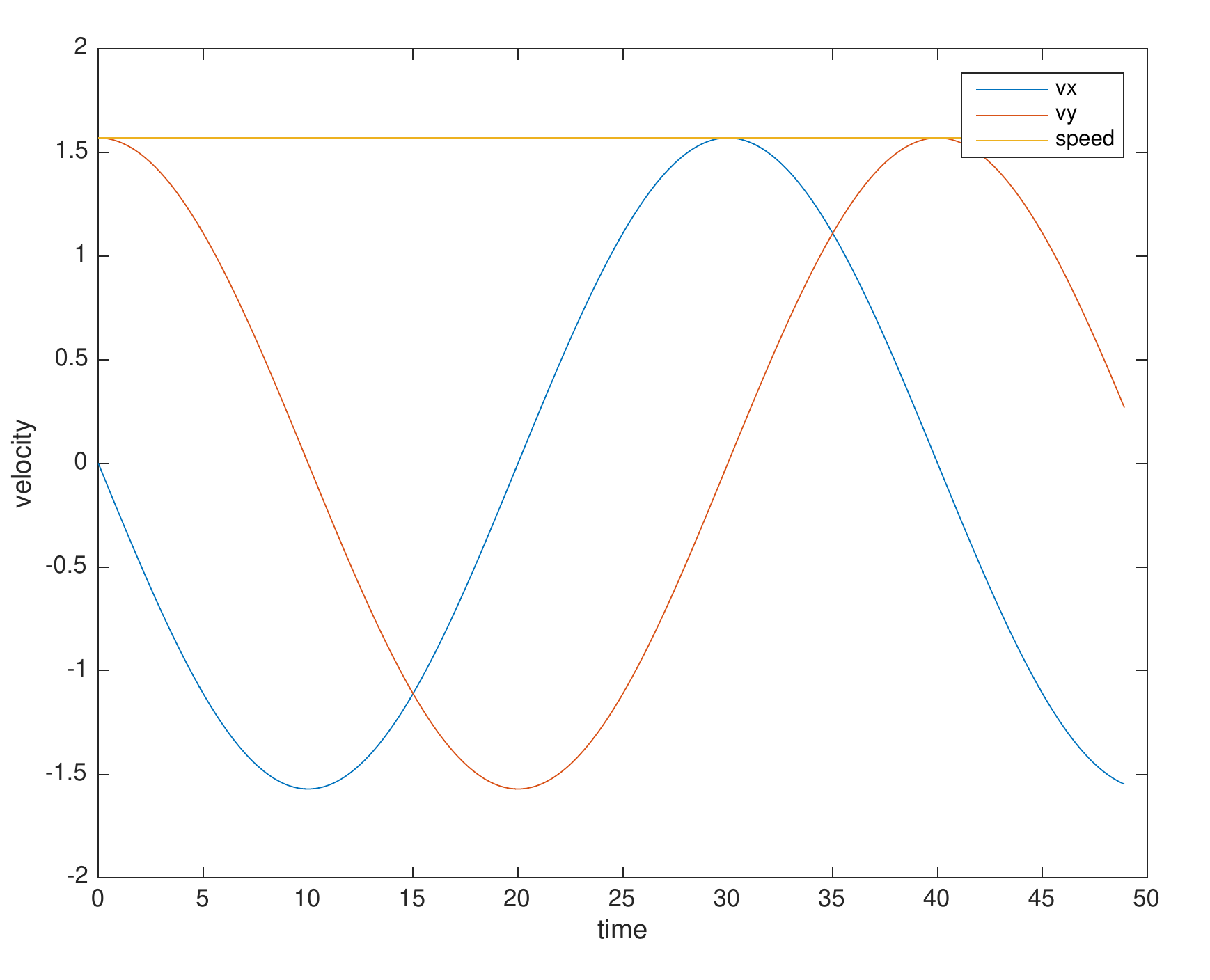}
    \caption{Velocity of circle trajectory}
\end{subfigure}
\caption{\label{circ_traj}Circle trajectory and velocity}
\end{figure}

The robot now travels along a circle of radius $10\,\rm{m}$, with constant speed that is approximately $1.57\,\rm{m/s}$. This is illustrated in Figure~\ref{circ_traj}. In contrast to the straight line case, we observed superior performance of the reduced EKF method over both the classical EKF and FE-EKF methods. 

\begin{figure}[H]
\centering
	\begin{subfigure}[b]{0.45\textwidth}
		\includegraphics[scale=0.4]{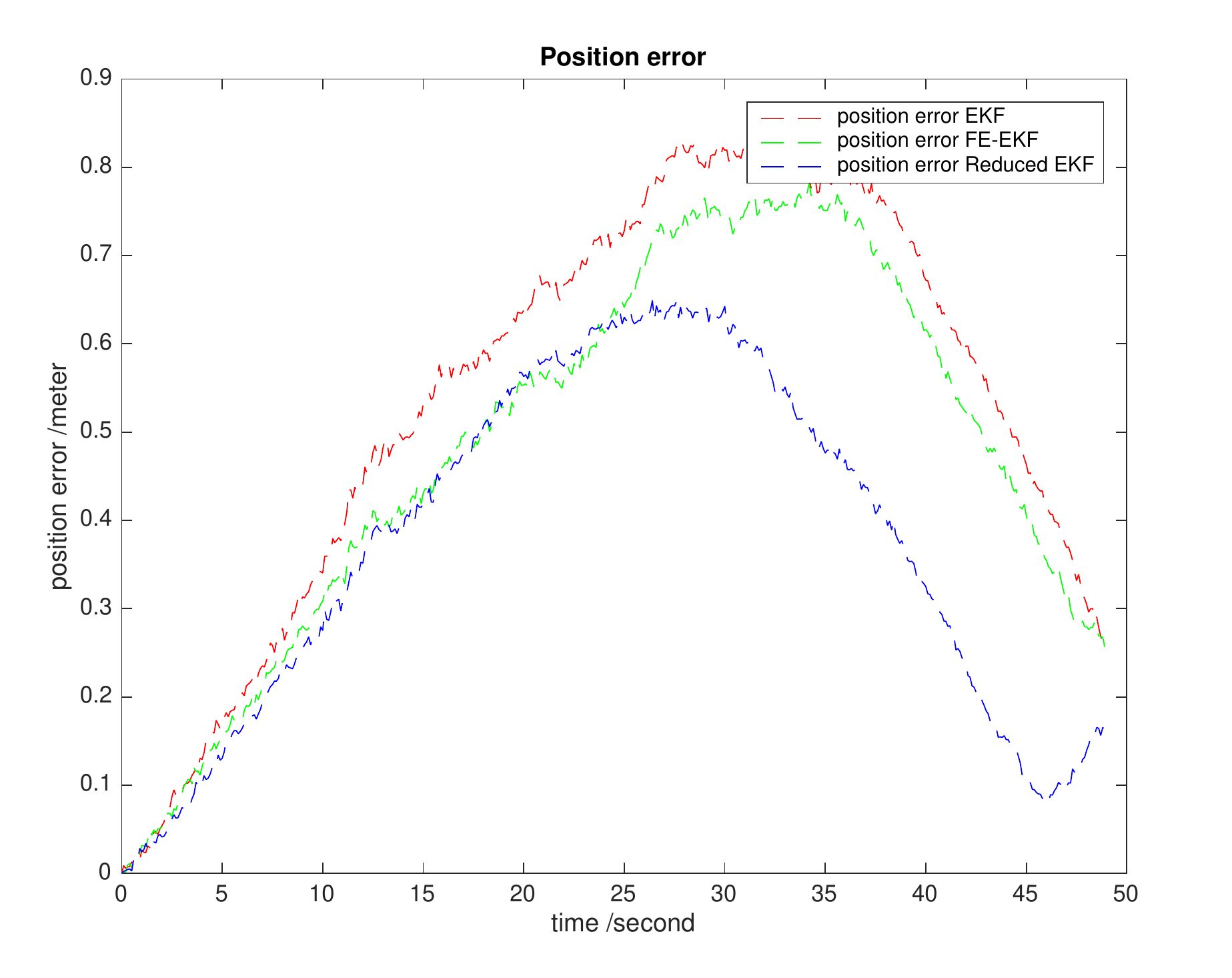}
		\caption{Position error, circle trajectory}
	\end{subfigure}
	\begin{subfigure}[b]{0.45\textwidth}
		\includegraphics[scale=0.4]{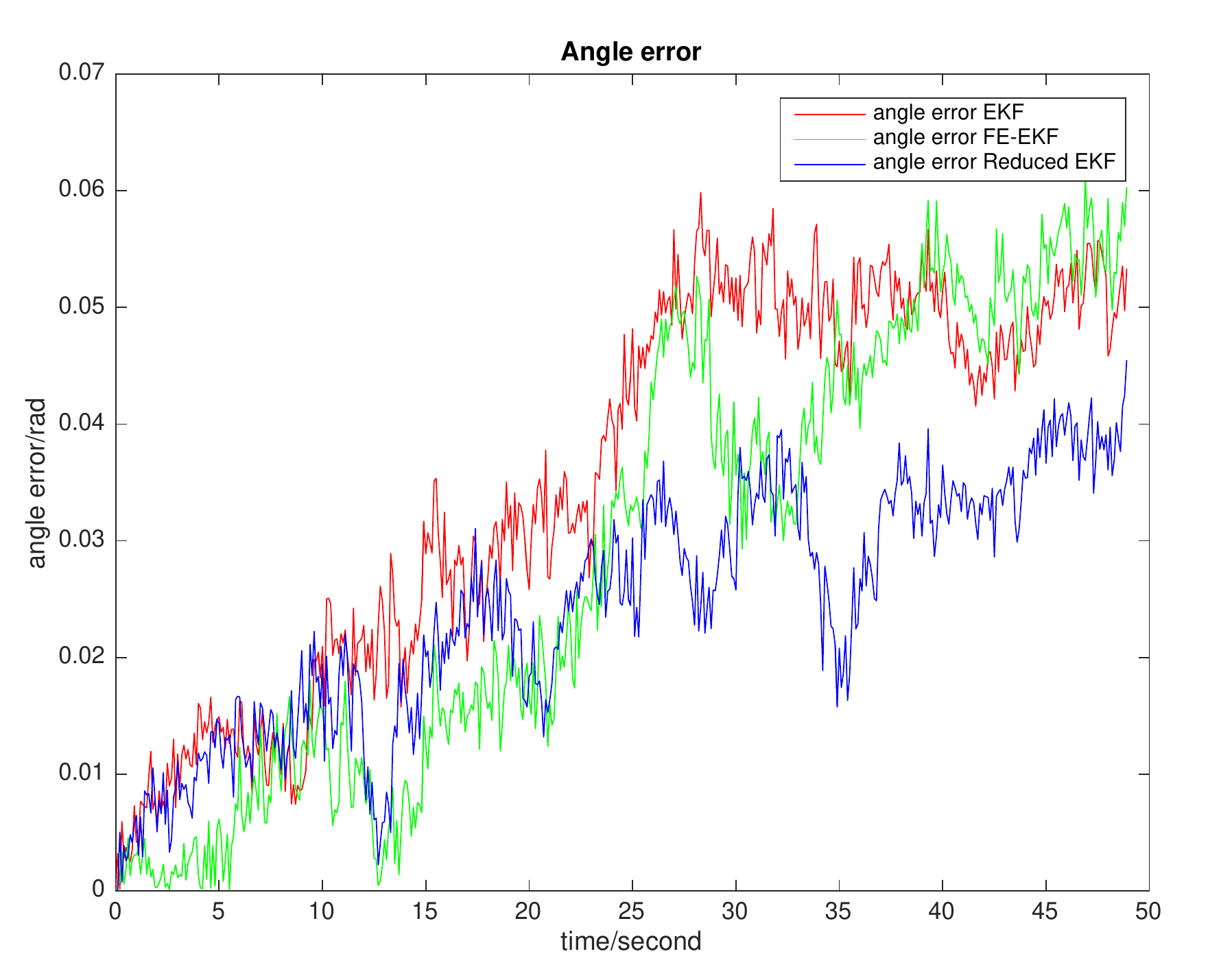}
		\caption{Angle error, circle trajectory}
\end{subfigure}
\caption{Error growth, circle trajectory\label{error_circle}}
\end{figure}

For the circle trajectory with low feature density, the average feature number per step is approximately 20, and for the high feature density, the average feature number per step is approximately 40. The average position error and orientation error of all three methods are summarized in Table~\ref{circle20} and Table~\ref{fig:circle40}. The improvement of reduced EKF method is quite significant compared to the classical EKF and FE-EKF methods. It can reduce the orientation error by approximately $50\%$ in each case, and as a consequence, the estimation of position is also improved, and the position error is reduced by approximately $30\%$ to $50\%$. A typical error growth for the circle trajectory is given in Figure~\ref{error_circle}.
\begin{table}\footnotesize
\centering
\begin{tabular}{|c|c|c|c|c|c|c|c|c|c|c|c|c|}
\hline
&  \multicolumn{12}{|c|}{Average error, circle trajectory, average feature num = 20} \\
\hline Frequency & \multicolumn{6}{|c|}{update frequency = 10Hz} & \multicolumn{6}{|c|}{update frequency = 20Hz} \\
\hline Method & \multicolumn{2}{|c|}{EKF} & \multicolumn{2}{|c|}{FE-EKF} & \multicolumn{2}{|c|}{Reduced EKF} & \multicolumn{2}{|c|}{EKF} & \multicolumn{2}{|c|}{FE-EKF} & \multicolumn{2}{|c|}{Reduced EKF} \\
\hline Error & $\delta{p}$ & $\delta{\theta}$ & $\delta{p}$ & $\delta{\theta}$ & $\delta{p}$ & $\delta{\theta}$ & $\delta{p}$ & $\delta{\theta}$ & $\delta{p}$ & $\delta{\theta}$ & $\delta{p}$ & $\delta{\theta}$ \\
\hline $Q_z=1\times 10^{-4}$ &  0.955 & 0.0875 & 1.063 & 0.0936 & 0.655 & 0.0443 & 0.746 & 0.0657 & 0.750 & 0.0576 & 0.532 & 0.0342\\
\hline $Q_z=2\times 10^{-4}$ &  1.630 & 0.152 & 1.496 & 0.119 & 0.959 & 0.0586 & 1.35 & 0.128 & 1.549 & 0.144 & 0.982 & 0.0667\\
\hline $Q_z=4\times 10^{-4}$ &  2.893 & 0.2914 & 2.721 & 0.2585 & 1.512 & 0.1033 & 2.295 & 0.2256 & 2.066 & 0.2102 & 1.469 & 0.1003\\
\hline
\end{tabular}
\caption{Average error, circle trajectory, average feature num = 20\label{circle20}}
\end{table}

\begin{table}\footnotesize
\centering
\begin{tabular}{|c|c|c|c|c|c|c|c|c|c|c|c|c|}
\hline
&  \multicolumn{12}{|c|}{Average error, circle trajectory, average feature num = 40} \\
\hline Frequency & \multicolumn{6}{|c|}{update frequency = 10Hz} & \multicolumn{6}{|c|}{update frequency = 20Hz} \\
\hline Method & \multicolumn{2}{|c|}{EKF} & \multicolumn{2}{|c|}{FE-EKF} & \multicolumn{2}{|c|}{Reduced EKF} & \multicolumn{2}{|c|}{EKF} & \multicolumn{2}{|c|}{FE-EKF} & \multicolumn{2}{|c|}{Reduced EKF} \\
\hline Error & $\delta{p}$ & $\delta{\theta}$ & $\delta{p}$ & $\delta{\theta}$ & $\delta{p}$ & $\delta{\theta}$ & $\delta{p}$ & $\delta{\theta}$ & $\delta{p}$ & $\delta{\theta}$ & $\delta{p}$ & $\delta{\theta}$ \\
\hline $Q_z=1\times 10^{-4}$ &  0.8161 & 0.0643 & 0.7972 & 0.0627 & 0.5815 & 0.0359 & 0.5843 & 0.0462 & 0.681 & 0.0530 & 0.568 & 0.0374\\
\hline $Q_z=2\times 10^{-4}$ &  1.2762 & 0.1193 & 1.3061 & 0.1091 & 0.9403 & 0.0591 & 1.2603 & 0.1164 & 1.2402 & 0.1039 & 0.9178 & 0.0573\\
\hline $Q_z=4\times 10^{-4}$ &  2.593 & 0.253 & 2.428 & 0.236 & 1.752 & 0.121 & 2.4 & 0.239 & 2.397 & 0.216 & 1.631 & 0.113\\
\hline
\end{tabular}
\caption{Average error, circle trajectory, average feature num = 40\label{fig:circle40}}
\end{table}

\subsection{General trajectory}
\begin{figure}[H]
\centering
	\begin{subfigure}[b]{0.45\textwidth}
		\includegraphics[scale=0.4]{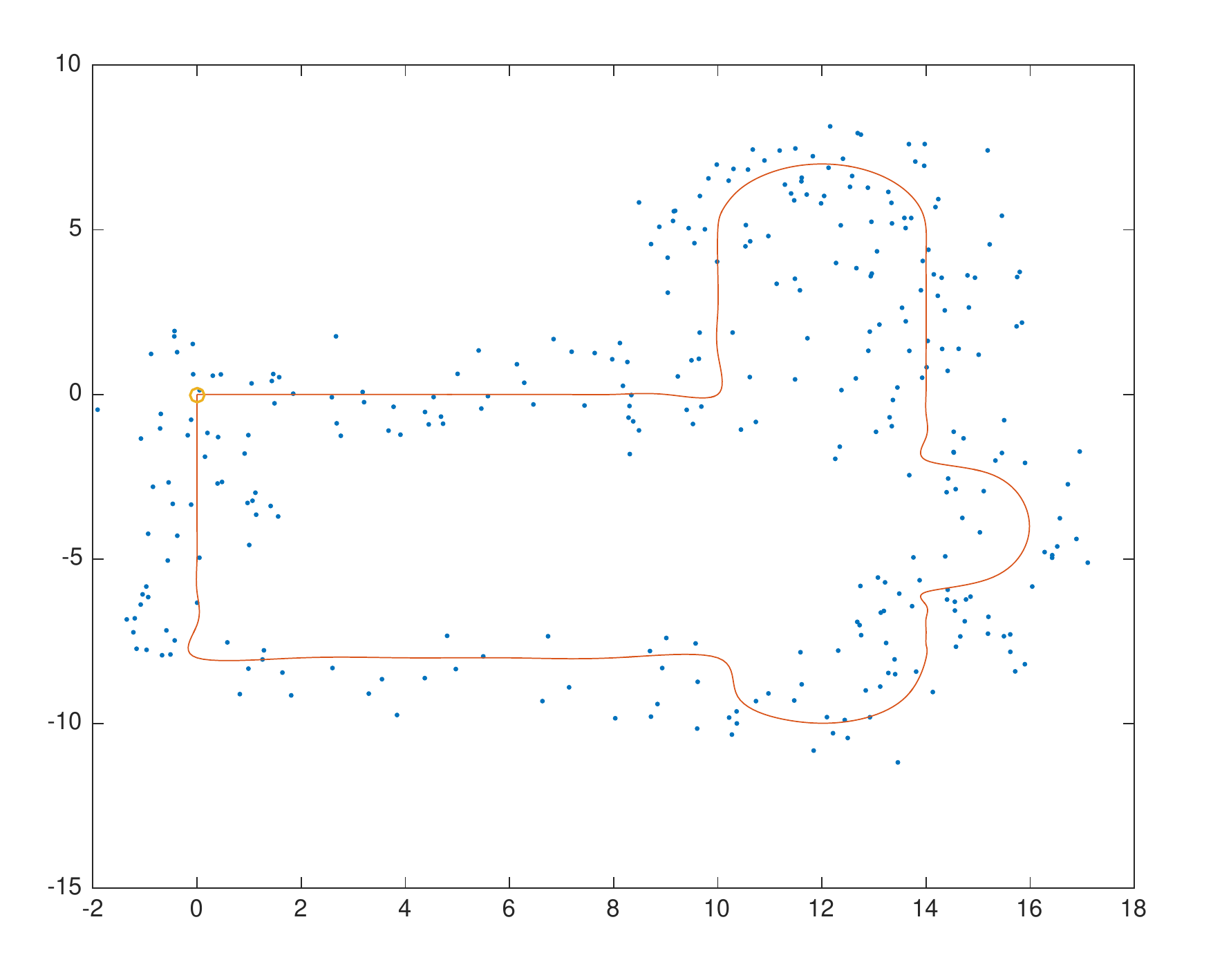}
		\caption{General, low density feature distribution}
	\end{subfigure}
	\begin{subfigure}[b]{0.45\textwidth}
		\includegraphics[scale=0.4]{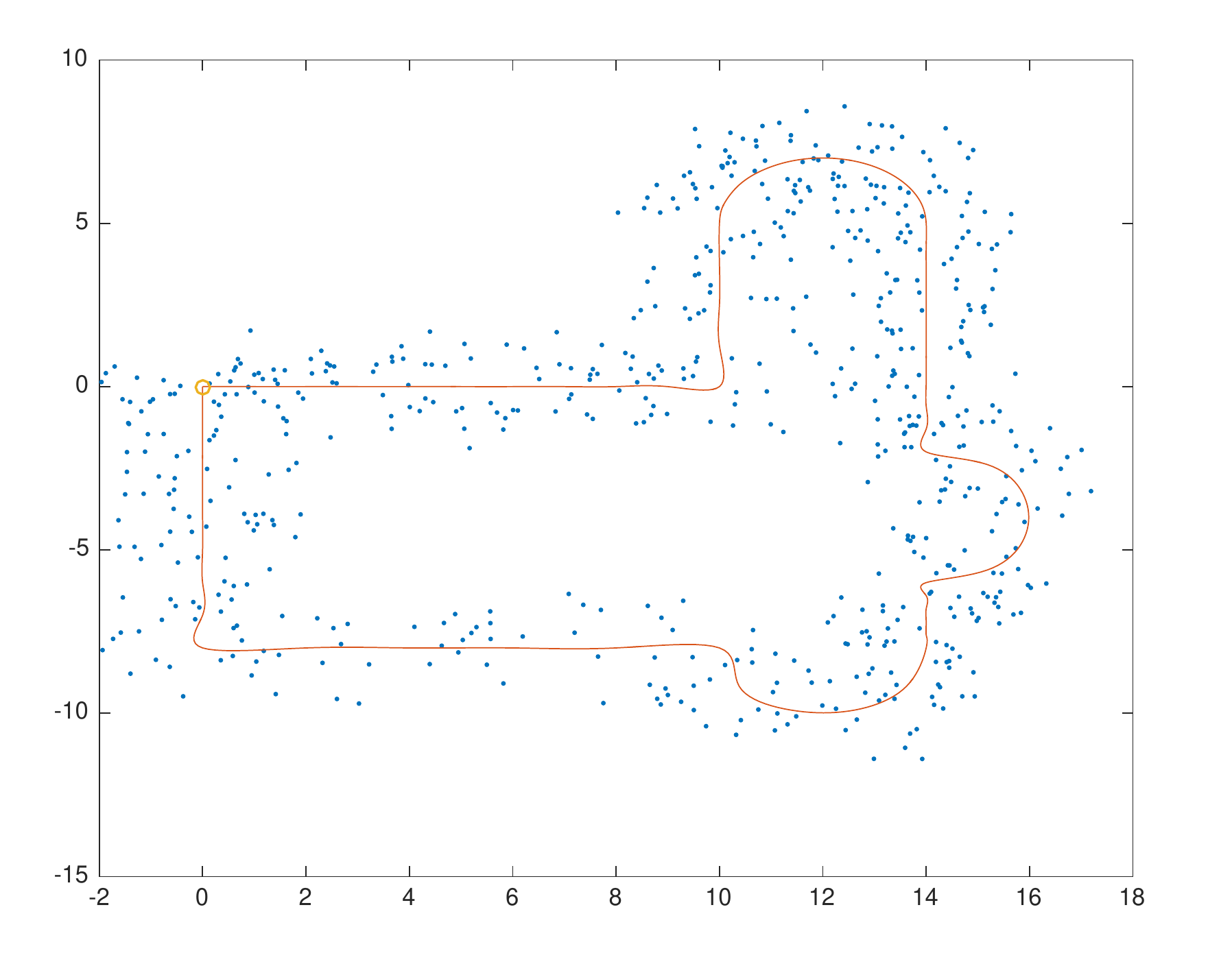}
		\caption{General, high density feature distribution}
\end{subfigure}
\begin{subfigure}[b]{0.45\textwidth}
  \includegraphics[scale=0.4]{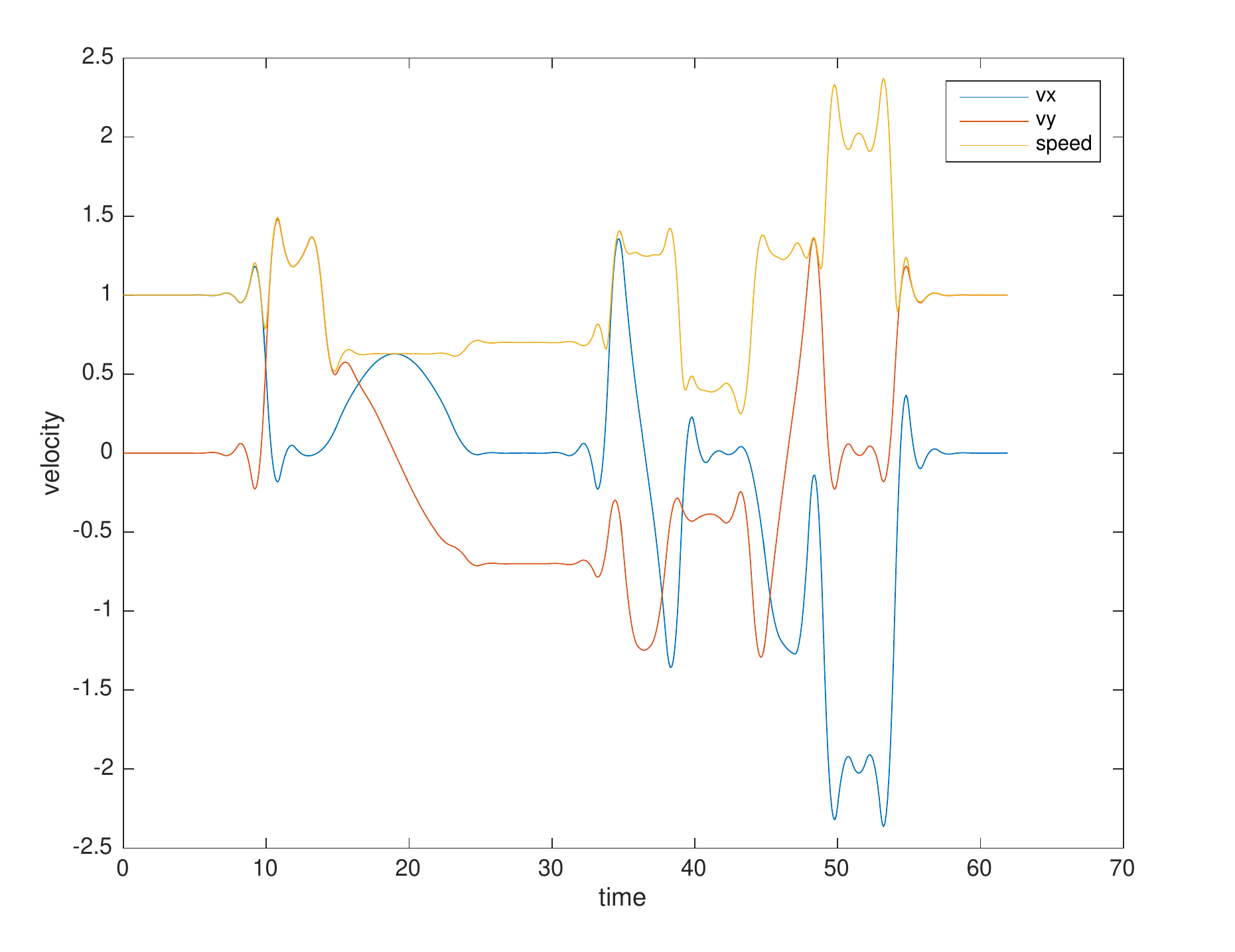}
    \caption{Velocity of general trajectory}
\end{subfigure} 
\caption{General trajectory and velocity\label{general_trajectory}}
\end{figure}

Besides the straight line and circle trajectories, we also consider a general trajectory, see Figure~\ref{general_trajectory}. The average feature number for the low feature density per step is approximately 25, and the average feature number for the high feature density per step is approximately 50.
As before, we tested all three methods on this general trajectory, and summarized their average position error and orientation error in Table~\ref{general25} and Table~\ref{general50}. When the update frequency is $20\,\rm{Hz}$, we observe noticeable improvement in both the position error and orientation error by the reduced EKF method. The FE-EKF method exhibits better performance than the classical EKF, but the reduced EKF method works even better. When the update frequence is $10\,\rm{Hz}$, and the measurement noise $Q_z=1\times 10^{-4}$, we see that the error of the reduced EKF is actually larger, which is not observed in the circle case. This might possibly be due to the high velocity of the general trajectory, since we construct the trajectory by cubic spline interpolation, and the velocity is very high around corners, which is unrealistic in practice as the actuation bounds in a robot result in bounded velocities as well.

\begin{figure}[H]
\centering
	\begin{subfigure}[b]{0.45\textwidth}
		\includegraphics[scale=0.4]{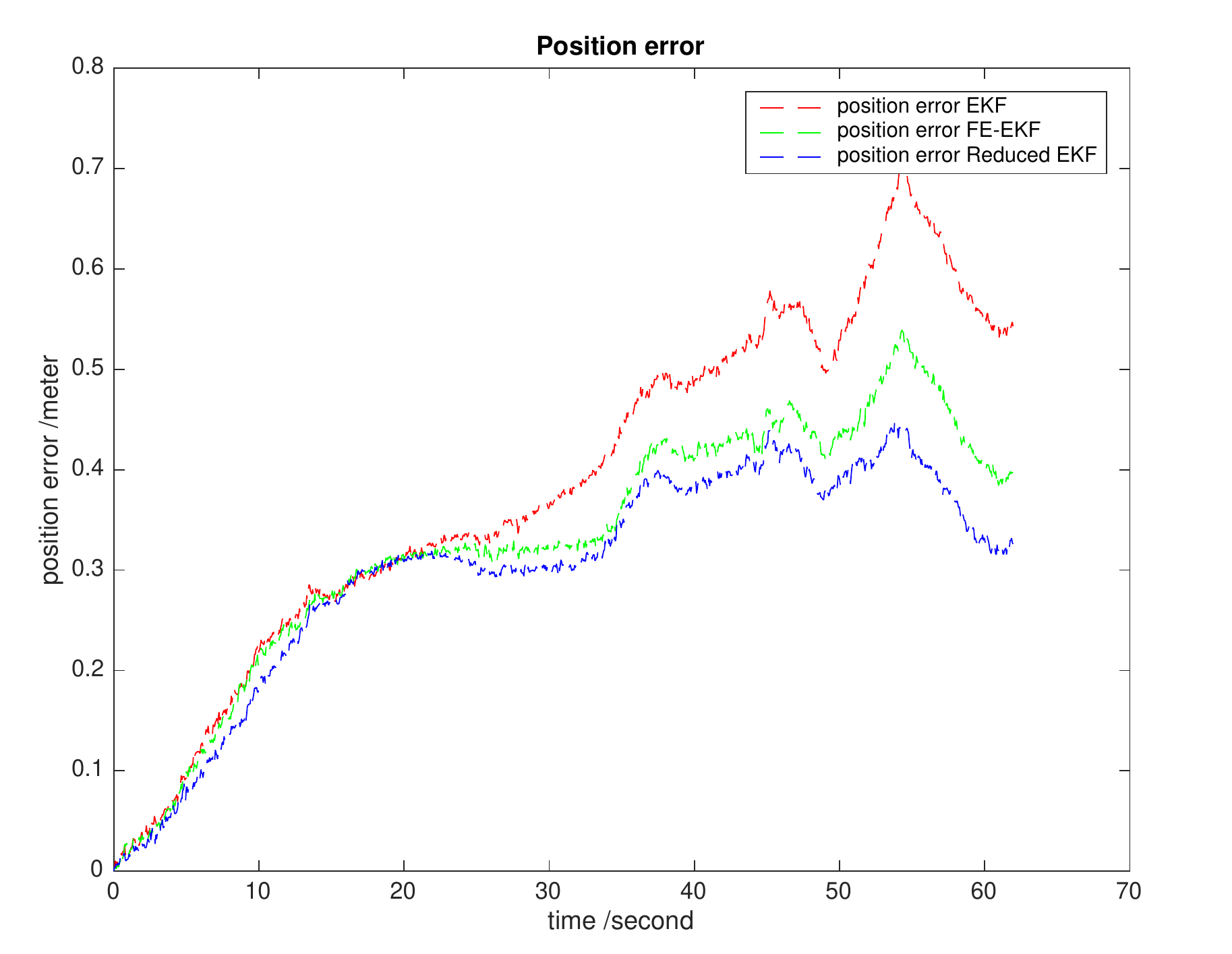}
		\caption{Position error, general trajectory}
	\end{subfigure}
	\begin{subfigure}[b]{0.45\textwidth}
		\includegraphics[scale=0.4]{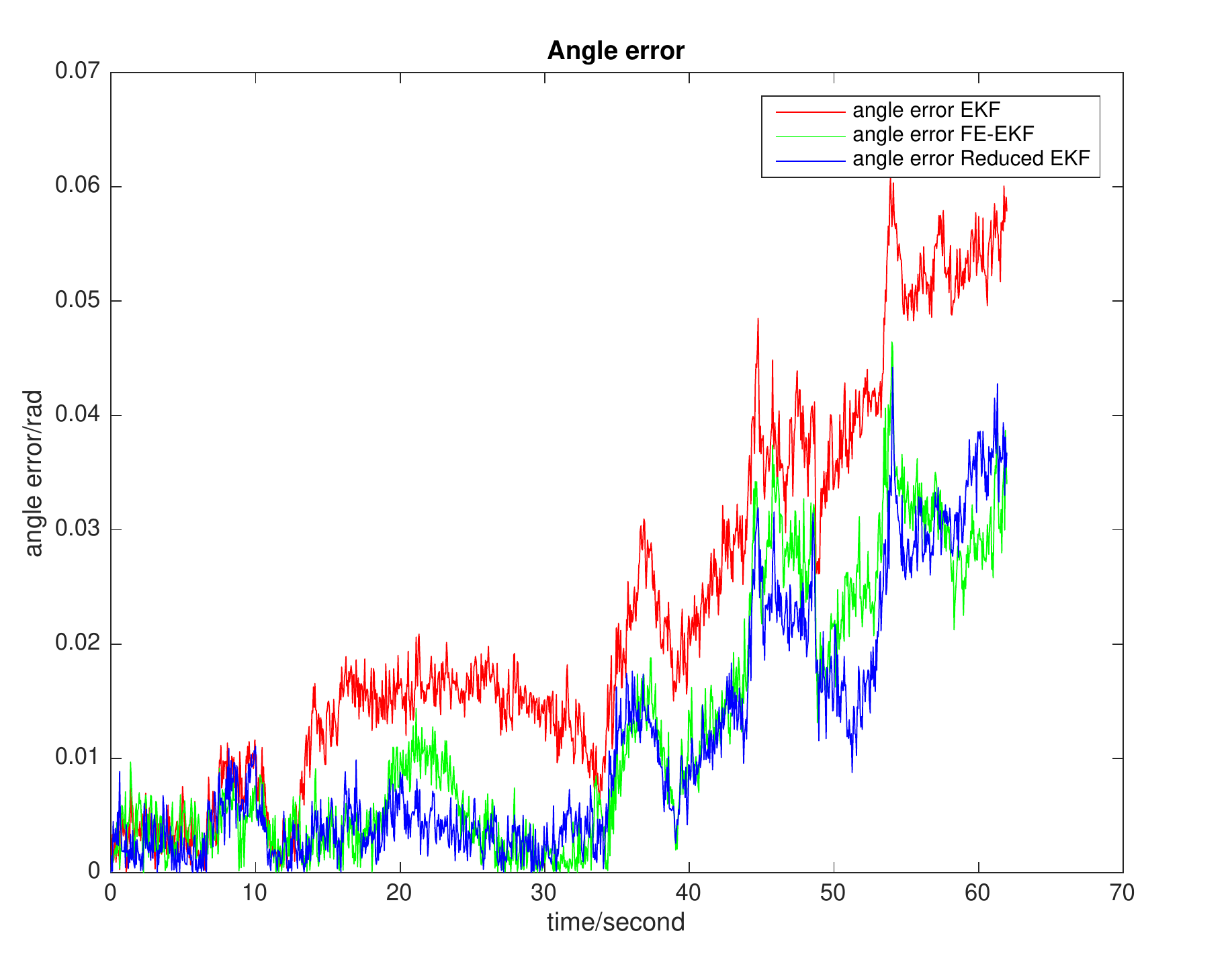}
		\caption{Angle error, general trajectory}
\end{subfigure}
\caption{Error growth, general trajectory\label{error_general}}
\end{figure}

Thus, the propagation error is large when the update frequency is low, and it seems that when the discrete system deviates far away from the true trajectory, forcing it to obey the observability property of the underlying continuous system actually degrades the performance of the algorithm. This is perhaps understandable, as when the timestep is too large to accurately integrate the trajectory, the unobservable subspace becomes poorly approximated as well. This is easily resolved by doubling the update frequency. Also, we see that even when the update frequency is $10\,\rm{Hz}$, when we increase the measurement noise to $4\times 10^{-4}$, the reduced EKF method maintains a rather robust orientation error, and its position error is now smaller than the classical EKF and FE-EKF methods, which demonstrates the stability of the reduced EKF method to large measurement noise. This can be understood in terms of the relative error from measurement noise compared to the propagation step. When the measurement noise is larger, the spurious information introduced to the unobservable variables will be significant, which degrades the estimation accuracy. We believe that the reduced EKF method that respects the observability constraint is best suitable for applications where the propagation step is more precise than the measurement step. This is the case in typical VIO systems, where the IMU for propagation has very high update frequency ($200$--$500\,\rm{Hz}$), and is very accurate for short time estimation, and the camera for measurement has low update frequency ($10$--$30\,\rm{Hz}$). A typical error growth in the general trajectory case is given in Figure~\ref{error_general}.

\begin{table}\footnotesize
\centering
\begin{tabular}{|c|c|c|c|c|c|c|c|c|c|c|c|c|}
\hline
&  \multicolumn{12}{|c|}{Average error, general trajectory, average feature num = 25} \\
\hline Frequency & \multicolumn{6}{|c|}{update frequency = 10Hz} & \multicolumn{6}{|c|}{update frequency = 20Hz} \\
\hline Method & \multicolumn{2}{|c|}{EKF} & \multicolumn{2}{|c|}{FE-EKF} & \multicolumn{2}{|c|}{Reduced EKF} & \multicolumn{2}{|c|}{EKF} & \multicolumn{2}{|c|}{FE-EKF} & \multicolumn{2}{|c|}{Reduced EKF} \\
\hline Error & $\delta{p}$ & $\delta{\theta}$ & $\delta{p}$ & $\delta{\theta}$ & $\delta{p}$ & $\delta{\theta}$ & $\delta{p}$ & $\delta{\theta}$ & $\delta{p}$ & $\delta{\theta}$ & $\delta{p}$ & $\delta{\theta}$ \\
\hline $Q_z=1\times 10^{-4}$ &  0.673 & 0.0789 & 0.635 & 0.0803 & 1.309 & 0.1062 & 0.6527 & 0.0698 & 0.5620 & 0.0605 & 0.4659 & 0.0371\\
\hline $Q_z=2\times 10^{-4}$ &  1.388 & 0.139 & 1.143 & 0.12 & 1.23 & 0.107 & 1.2575 & 0.1156 & 0.9985 & 0.0997 & 0.8022 & 0.0654\\
\hline $Q_z=4\times 10^{-4}$ &  2.459 & 0.227 & 2.111 & 0.2267 & 1.592 & 0.1174 & 2.3852 & 0.2039 & 1.8924 & 0.1794 & 1.4595 & 0.1042\\
\hline
\end{tabular}
\caption{Average error, general trajectory, average feature num = 25\label{general25}}
\end{table}

\begin{table}\footnotesize
\centering
\begin{tabular}{|c|c|c|c|c|c|c|c|c|c|c|c|c|}
\hline
&  \multicolumn{12}{|c|}{Average error, general trajectory, average feature num = 50} \\
\hline Frequency & \multicolumn{6}{|c|}{update frequency = 10Hz} & \multicolumn{6}{|c|}{update frequency = 20Hz} \\
\hline Method & \multicolumn{2}{|c|}{EKF} & \multicolumn{2}{|c|}{FE-EKF} & \multicolumn{2}{|c|}{Reduced EKF} & \multicolumn{2}{|c|}{EKF} & \multicolumn{2}{|c|}{FE-EKF} & \multicolumn{2}{|c|}{Reduced EKF} \\
\hline Error & $\delta{p}$ & $\delta{\theta}$ & $\delta{p}$ & $\delta{\theta}$ & $\delta{p}$ & $\delta{\theta}$ & $\delta{p}$ & $\delta{\theta}$ & $\delta{p}$ & $\delta{\theta}$ & $\delta{p}$ & $\delta{\theta}$ \\
\hline $Q_z=1\times 10^{-4}$ &  0.6039 & 0.0609 & 0.5815 & 0.0636 & 1.308 & 0.1049 & 0.5518 & 0.0518 & 0.4810 & 0.0447 & 0.4538 & 0.0381\\
\hline $Q_z=2\times 10^{-4}$ &  1.2201 & 0.1047 & 0.9836 & 0.0945 & 1.0738 & 0.0838 & 1.1676 & 0.1021 & 0.9752 & 0.0869 & 0.73 & 0.0611\\
\hline $Q_z=4\times 10^{-4}$ &  2.4264 & 0.2126 & 2.0421 & 0.1904 & 1.4674 & 0.1121 & 2.447 & 0.2141 & 1.8347 & 0.145 & 1.4322 & 0.1165\\
\hline
\end{tabular}
\caption{Average error, general trajectory, average feature num = 50\label{general50}}
\end{table}

\section{Conclusion and Future Research}
We consider the observability properties of a control system in the context of Lie group symmetries. Due to the Lie group symmetry, the state variable can be decomposed into the unobservable and observable parts $(x_N, x_O)$ explicitly.  Since the Lie group action on $x_N$ is trivial, the problem of requiring that the propagation step of the Kalman filter satisfies the geometric property that the unobservable subspace remains perpendicular to the measurement $dh$ is naturally formulated in terms of the reduced control system on the symmetry reduced space corresponding to the observable variables.  Moreover, in order to deal with the problem that this property is destroyed at the measurement update step of the Kalman filter, we propose a reduced Bayesian inference method, i.e., only the observable part $x_O$ is updated by the reduced measurement. This special procedure guarantees that the unobservable variable $x_N$ remains unobservable, and no spurious information is introduced by the measurement update.  We applied this approach to the problem of the planar robot equipped with odometry sensors and cameras, and the reduced Kalman filter method based on observability considerations outperforms the traditional Kalman filter as well as the FEJ-Kalman filter by quite a lot.

 In the future, we intend to apply this technique to a practical VIO system, where the geometry involved is more complicated. In particular, the rotation group in $\mathbb{R}^3$ is $SO(3)$ which is nonabelian, and only rotations about the gravity direction are unobservable. Also, we would like to update the multi-state constraint Kalman filter (MSCKF) framework, where the state variables include historic poses of the camera, which makes the decomposition into observable and unobservable components more interesting. In principle, all traditional techniques that are applicable to the Kalman filter, such as UKF, particle filters, and so on, could all be combined with our reduced update method based on the observability property, and we would like to test the efficiency of our method for all related practical applications.

\section*{Acknowledgements}
This research has been supported in part by NSF under grants DMS-1010687, CMMI-1029445, DMS-1065972, CMMI-1334759, DMS-1411792, DMS-1345013, DMS-1813635, and by AFOSR under grant FA9550-18-1-0288.

\bibliography{reference}
\bibliographystyle{plainnat}
\end{document}